\documentclass{amsart}

\usepackage{amssymb}
\usepackage{amsmath}
\usepackage{amsthm}
\usepackage{amsfonts}

\usepackage{a4wide}
\usepackage{graphicx}
\usepackage[usenames,dvipsnames]{xcolor}

\newtheorem{theorem}{Theorem}

\newtheorem{corollary}[theorem]{Corollary}

\theoremstyle{definition}
\newtheorem{definition}[theorem]{Definition}

\newtheorem{example}[theorem]{Example}

\begin{document}


\title
[Finite non-metabelian Schur \(\sigma\)-Galois groups]
{Finite non-metabelian Schur \(\sigma\)-Galois groups \\ of class field towers}

\author{Daniel C. Mayer}
\address{Naglergasse 53 \\ 8010 Graz \\ Austria}
\email{algebraic.number.theory@algebra.at}
\urladdr{http://www.algebra.at}

\thanks{Research supported by the Austrian Science Fund (FWF): P 26008-N25 and J0497-PHY, and by the EUREA}


\date{October 21, 2020}


\begin{abstract}
For each odd prime \(p\ge 5\),
there exist finite \(p\)-groups \(G\)
with derived quotient \(G/G^\prime\simeq C_p\times C_p\)
and nearly constant transfer kernel type \(\varkappa(G)=(\varkappa_1,\ldots,\varkappa_{p+1})\)
having two fixed points \(\varkappa_1=1\), \(\varkappa_2=2\), and \(\varkappa_i=2\) for \(3\le i\le p+1\).
It is proved that, for \(p=7\), this type \(\varkappa(G)\)
with the simplest possible case of logarithmic abelian quotient invariants
\(\tau(G)=(H_i/H_i^\prime)_{1\le i\le 8}=(1^5,1^3,(21)^6)\)
of the eight maximal subgroups \(H_i\)
is realized by exactly \(98\) non-metabelian Schur \(\sigma\)-groups \(S\) of order \(7^{11}\)
with fixed derived length \(\mathrm{dl}(S)=3\)
and metabelianizations \(S/S^{\prime\prime}\) of order \(7^7\).
For \(p=5\), the type \(\varkappa(G)\) with
\(\tau(G)=(H_i/H_i^\prime)_{1\le i\le 6}=(21^3,1^3,(21)^4)\)
leads to infinitely many non-metabelian Schur \(\sigma\)-groups \(S\) of order at least \(5^{14}\)
with unbounded derived length \(\mathrm{dl}(S)\ge 3\)
and metabelianizations \(S/S^{\prime\prime}\) of fixed order \(5^7\).
These results admit the conclusion that
\(d=-159\,592\) is the first known discriminant of an
imaginary quadratic field with \(7\)-class field tower
of precise length \(\ell_7(K)=3\),
and \(d=-90\,868\) is a discriminant of an
imaginary quadratic field with \(5\)-class field tower
of length \(\ell_5(K)\ge 3\), whose exact length remains unknown.
\end{abstract}

\maketitle



\section{Introduction}
\label{s:Intro}

\noindent
Let \(p\) be an odd prime number.
Using the Shafarevich theorem
\cite[Thm. 5.1, p. 28]{Ma2015d}
on the relation rank \(d_2=\dim_{\mathbb{F}_p}(\mathrm{H}^2(G,\mathbb{F}_p))\)
of the automorphism group \(G=\mathrm{Gal}(\mathrm{F}_p^\infty(K)/K)\)
of the maximal unramified pro-\(p\)-extension \(\mathrm{F}_p^\infty(K)\),
which is called the \textit{Hilbert \(p\)-class field tower},
of an algebraic number field \(K\),
Koch and Venkov
\cite{KoVe}
showed that for an imaginary quadratic field \(K=\mathbb{Q}(\sqrt{d})\)
with fundamental discriminant \(d<0\)
the Galois group \(G\) must be a Schur \(\sigma\)-group
\cite{Ag},
and that this pro-\(p\)-group can only be a finite \(p\)-group
when the \(p\)-class rank \(\varrho_p=\dim_{\mathbb{F}_p}(\mathrm{Cl}_p(K)/\mathrm{Cl}_p(K)^p)\) of \(K\)
is bounded by \(\varrho_p\le 2\).

The aim of the present work is to grant access to the most recent information on
finite non-metabelian Schur \(\sigma\)-groups \(S\) with odd prime power order \(\#S=p^n\),
and their actual realization \(S\simeq\mathrm{Gal}(\mathrm{F}_p^\infty(K)/K)\)
as Galois groups of finite \(p\)-class field towers with at least three stages
over imaginary quadratic fields \(K\).
In view of the necessary condition \(\varrho_p\le 2\)
we restrict our investigations to elementary bicyclic \(p\)-class groups \(\mathrm{Cl}_p(K)\simeq C_p\times C_p\).

In order to reduce the probability of infinite families
\cite{BaBu}
of Schur \(\sigma\)-groups \(S\) with unbounded derived length \(\mathrm{dl}(S)\),
we recall the following terminology
\cite{Ma2015a,Ma2015c,Ma2016a,Ma2018b,Ma2018c}. 


\begin{definition}
\label{dfn:Periodic}
A finite \(p\)-group is called \textbf{periodic}
if its isomorphism class is vertex of an infinite descendant tree with fixed coclass.
Otherwise it is called \textbf{sporadic}.
\end{definition}


\noindent
Unbounded descendant trees
whose edges are restricted to step size \(s=1\),
and whose vertices consequently share a common coclass,
are called \textit{coclass trees}
\cite{Ma2015a}.
We shall focus our attention on Schur \(\sigma\)-groups \(S\)
with periodic metabelianizations \(M=S/S^{\prime\prime}\) on coclass trees.

\newpage


\section{Theoretical foundations}
\label{s:Foundations}

\noindent
For the convenience of the reader we begin with a collection of concepts
which appear in the title of this work: non-metabelian groups, Schur-groups, and \(\sigma\)-groups.


The \textit{derived series} \(G^{(0)}\ge G^{(1)}\ge G^{(2)}\ge\ldots\ge G^{(n)}\ge\ldots \) of a group \(G\)
is defined recursively by
\(G^{(0)}:=G\) and \((\forall\,i\in\mathbb{N})\,G^{(i)}:=\lbrack G^{(i-1)},G^{(i-1)}\rbrack\).
If this series becomes stationary in the trivial group \(1\),
that is, either \(G^{(0)}=1\) or \((\exists\,n\in\mathbb{N})\,G^{(n-1)}>G^{(n)}=1\),
then \(G\) is called \textit{solvable}.
The \textit{derived length} of \(G\) is defined by
\(\mathrm{dl}(G):=\inf\lbrace i\in\mathbb{N}_0\mid G^{(i)}=1\rbrace\in\mathbb{N}_0\cup\lbrace\infty\rbrace\).
This invariant admits a coarse classification of groups:

\begin{itemize}
\item
\(\mathrm{dl}(G)=0\) \(\Longleftrightarrow\) \(G=1\) is the \textit{trivial} group.
\item
\(\mathrm{dl}(G)=1\) \(\Longleftrightarrow\) \(G>G^\prime=1\), i.e. \(G\) is (non-trivial) \textit{abelian}.
\item
\(\mathrm{dl}(G)=2\) \(\Longleftrightarrow\) \(G>G^\prime>G^{\prime\prime}=1\), i.e. \(G^\prime\) is abelian, and \(G\) is non-abelian but \textit{metabelian}.
\item
\(\mathrm{dl}(G)=\infty\) \(\Longleftrightarrow\) \((\forall\,i\in\mathbb{N}_0)\,G^{(i)}>1\)
 \(\Longleftrightarrow\) \(G\) is \textit{non-solvable}.
\end{itemize}

\begin{definition}
\label{dfn:NonMetabelian}
A finite group \(G\) is called \textit{non-metabelian}
if it possesses derived length \(\mathrm{dl}(G)\ge 3\).
\end{definition}


Let \(p\) be a prime number, and \(G\) be a pro-\(p\) group.
Then the \textit{generator rank} \(d_1(G):=\dim_{\mathbb{F}_p}(\mathrm{H}^1(G,\mathbb{F}_p))\)
and the \textit{relation rank} \(d_2(G):=\dim_{\mathbb{F}_p}(\mathrm{H}^2(G,\mathbb{F}_p))\) 
are two cohomological invariants of \(G\).

\begin{definition}
\label{dfn:Schur}
A pro-\(p\) group \(G\) is called a \textit{Schur-group}
if \(d_2(G)=d_1(G)\).
In this case \(G\) is said to possess a \textit{balanced presentation},
since the relation rank equals the generator rank
\cite{Ag,KoVe}.
\end{definition}


Let \(\sigma\in\mathrm{Aut}(G)\) be an automorphism of a group \(G\).
Since the derived subgroup \(G^\prime\) of \(G\) is characteristic,
\(\sigma\) induces an automorphism on the abelianization \(G/G^\prime\) of \(G\).
If the induced automorphism is denoted by \(\hat{\sigma}\), then
\((\forall\,n\in\mathbb{N})\,(\forall\,x\in G)\,\hat{\sigma}^n(xG^\prime)=\sigma^n(x)G^\prime\),
which shows that the order of \(\hat{\sigma}\) divides the order of \(\sigma\),
if the latter is finite.

\begin{definition}
\label{dfn:Sigma}
A group \(G\) is called a \textit{\(\sigma\)-group}
if it possesses an automorphism \(\sigma\in\mathrm{Aut}(G)\)
which acts as inversion on the abelianization, that is,
\((\forall\,x\in G)\,\sigma(x)G^\prime=x^{-1}G^\prime\).
\end{definition}

Observe that a group \(G\) is abelian if and only if
the inversion is an automorphism of \(G\).
Consequently,
every abelian group is a \(\sigma\)-group, and
a finite \(\sigma\)-group \(G\) has an automorphism group with even order \(\#\mathrm{Aut}(G)\),
since the inversion on \(G/G^\prime\), induced by \(\sigma\), is an involution of order \(2\).



\section{\(7\)-class field towers with precise length three}
\label{s:7ClassTowers}

\noindent
In this section we present our main result
concerning the discovery of the first three-stage towers of \(7\)-class fields,
which do not occur in the previous literature.
In order to ensure the exact length three,
we abstain from transfer kernel types (TKT) of sporadic groups,
which are permutations.


\begin{theorem}
\label{thm:7ClassTowers} 
(Sufficient criterion for a \(7\)-class field tower with exactly three stages) \\
An imaginary quadratic number field \(K=\mathbb{Q}(\sqrt{d})\)
with fundamental discriminant \(d<0\),
elementary bicyclic \(7\)-class group \(\mathrm{Cl}_7(K)\simeq C_7\times C_7\),
and Artin pattern \(\mathrm{AP}(K)=(\varkappa(K),\tau(K))\)
with capitulation type \(\varkappa(K)\sim (1,2,2,2,2,2,2,2)\)
in the eight unramified cyclic septic extensions \((E_i)_{1\le i\le 8}\) of \(K\)
and logarithmic abelian type invariants \(\tau(K)\sim (1^5,1^3,21,21,21,21,21,21)\)
of the \(7\)-class groups \((\mathrm{Cl}_7(E_i))_{1\le i\le 8}\)
possesses a finite \(7\)-class field tower
\(K<\mathrm{F}_7^1(K)<\mathrm{F}_7^2(K)<\mathrm{F}_7^3(K)=\mathrm{F}_7^\infty(K)\)
with precise length \(\ell_7(K)=3\).
\end{theorem}


\noindent
Now we supplement the purely arithmetical information of Theorem
\ref{thm:7ClassTowers}
with group theoretical data.
We use the notation of the SmallGroups database
\cite{BEO1,BEO2}
and the ANUPQ package
\cite{GNO}.


\begin{corollary}
\label{cor:7ClassTowers} 
Under the conditions of Theorem
\ref{thm:7ClassTowers},
there are only two possibilities for the metabelian second \(7\)-class group
\(M=\mathrm{Gal}(\mathrm{F}_7^2(K)/K)\) of \(K\), namely
\(M\simeq\mathrm{SmallGroup}(823\,543,m)\) with either \(m=1990\) or \(m=1991\).
Both share the invariants \(\mathrm{lo}=7\), \(\mathrm{cl}=5\), \(\mathrm{cc}=2\), \(d_2=3\).
See Figure
\ref{fig:SchurSigmaLogOrd11}.

There are \(98\) candidates for the non-metabelian Schur \(\sigma\)-Galois group 
\(S=\mathrm{Gal}(\mathrm{F}_7^\infty(K)/K)\)
of the \(7\)-class field tower of \(K\).
In dependence on the metabelianization \(M\simeq S/S^{\prime\prime}\),
they are given by

\begin{equation}
\label{eqn:7ClassTowers}
\begin{aligned}
S &\simeq\mathrm{SmallGroup}(117\,649,708)-\#2;\ell-\#1,k-\#2;j \text{ with } 1\le j\le 7 \text{ and} \\
\ell &\in\lbrace 4+17v\mid 0\le v\le 6\rbrace,\ k=\begin{cases}1 \text{ for } \ell=4, \\2 \text{ otherwise,}\end{cases} \text{ if } m=1990, \\
\ell &\in\lbrace 11+17v\mid 0\le v\le 6\rbrace,\ k=\begin{cases}1 \text{ for } \ell=11, \\4 \text{ otherwise,}\end{cases} \text{ if } m=1991. \\
\end{aligned}
\end{equation}

\noindent
All of them share the invariants \(\mathrm{lo}=11\), \(\mathrm{cl}=7\), \(\mathrm{cc}=4\),
and \(d_2=2=d_1=\dim_{\mathbb{F}_7}(\mathrm{H}^1(S,\mathbb{F}_7))=\varrho_7\).
\end{corollary}


\begin{proof}
The next section is devoted to the rigorous justification of Theorem
\ref{thm:7ClassTowers}
and Corollary
\ref{cor:7ClassTowers}.
\end{proof}


\begin{example}
\label{exm:7ClassTowers}
The three absolutely smallest discriminants of imaginary quadratic number fields \(K=\mathbb{Q}(\sqrt{d})\)
which satisfy the conditions in Theorem
\ref{thm:7ClassTowers}
are \(d=-159\,592\), \(d=-611\,076\), and \(d=-839\,147\).
The \(7\)-class field tower over these three fields has exactly three stages.
\end{example} 



\section{\(7\)-groups connected with \(7\)-class towers}
\label{s:7Groups}

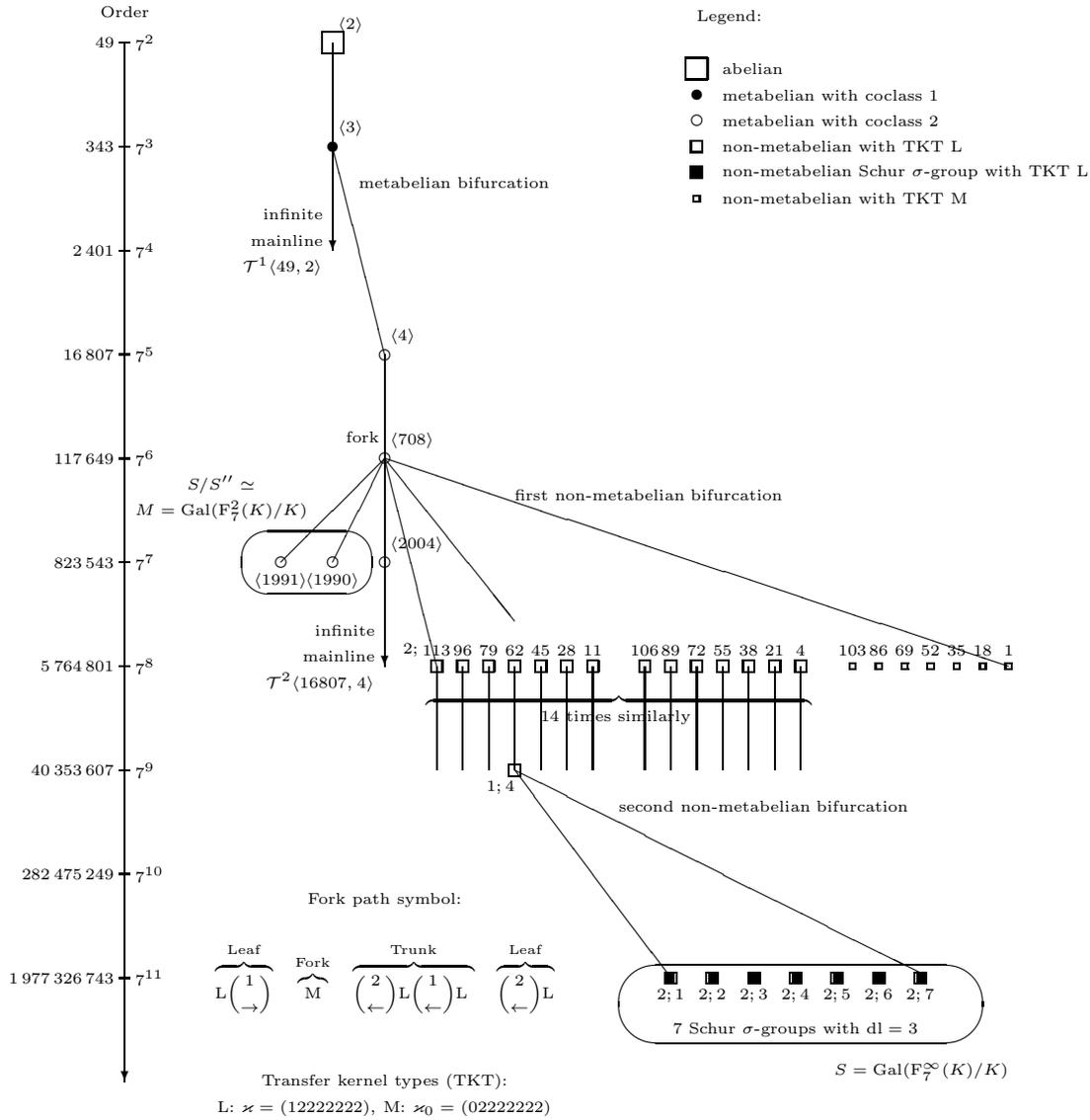
\begin{figure}[ht]
\caption{\(7\)-class field tower groups of order \(7^{11}\) on the tree \(\mathcal{T}\langle 16\,807,4\rangle\)}
\label{fig:SchurSigmaLogOrd11}


\setlength{\unitlength}{0.7cm}
\begin{picture}(18,21)(-6,-20)

{\tiny

\put(-5,0.5){\makebox(0,0)[cb]{Order}}
\put(-5,0){\line(0,-1){18}}
\multiput(-5.1,0)(0,-2){10}{\line(1,0){0.2}}
\put(-5.2,0){\makebox(0,0)[rc]{\(49\)}}
\put(-4.8,0){\makebox(0,0)[lc]{\(7^2\)}}
\put(-5.2,-2){\makebox(0,0)[rc]{\(343\)}}
\put(-4.8,-2){\makebox(0,0)[lc]{\(7^3\)}}
\put(-5.2,-4){\makebox(0,0)[rc]{\(2\,401\)}}
\put(-4.8,-4){\makebox(0,0)[lc]{\(7^4\)}}
\put(-5.2,-6){\makebox(0,0)[rc]{\(16\,807\)}}
\put(-4.8,-6){\makebox(0,0)[lc]{\(7^5\)}}
\put(-5.2,-8){\makebox(0,0)[rc]{\(117\,649\)}}
\put(-4.8,-8){\makebox(0,0)[lc]{\(7^6\)}}
\put(-5.2,-10){\makebox(0,0)[rc]{\(823\,543\)}}
\put(-4.8,-10){\makebox(0,0)[lc]{\(7^7\)}}
\put(-5.2,-12){\makebox(0,0)[rc]{\(5\,764\,801\)}}
\put(-4.8,-12){\makebox(0,0)[lc]{\(7^8\)}}
\put(-5.2,-14){\makebox(0,0)[rc]{\(40\,353\,607\)}}
\put(-4.8,-14){\makebox(0,0)[lc]{\(7^9\)}}
\put(-5.2,-16){\makebox(0,0)[rc]{\(282\,475\,249\)}}
\put(-4.8,-16){\makebox(0,0)[lc]{\(7^{10}\)}}
\put(-5.2,-18){\makebox(0,0)[rc]{\(1\,977\,326\,743\)}}
\put(-4.8,-18){\makebox(0,0)[lc]{\(7^{11}\)}}
\put(-5,-18){\vector(0,-1){2}}


\put(6,0.5){\makebox(0,0)[lc]{Legend:}}
\put(5.8,-0.7){\framebox(0.4,0.4){}}
\put(6,-1){\circle*{0.2}}
\put(6,-1.5){\circle{0.2}}
\put(5.9,-2.1){\framebox(0.2,0.2){}}
\put(5.9,-2.6){\framebox(0.2,0.2){\vrule height 3pt width 3pt}}
\put(5.95,-3.05){\framebox(0.1,0.1){}}

\put(6.5,-0.5){\makebox(0,0)[lc]{abelian}}
\put(6.5,-1){\makebox(0,0)[lc]{metabelian with coclass 1}}
\put(6.5,-1.5){\makebox(0,0)[lc]{metabelian with coclass 2}}
\put(6.5,-2){\makebox(0,0)[lc]{non-metabelian with TKT L}}
\put(6.5,-2.5){\makebox(0,0)[lc]{non-metabelian Schur \(\sigma\)-group with TKT L}}
\put(6.5,-3){\makebox(0,0)[lc]{non-metabelian with TKT M}}


\put(0,-16.5){\makebox(0,0)[cc]{Fork path symbol:}}
\put(0,-18){\makebox(0,0)[cc]{
\(\overbrace{\mathrm{L}\binom{1}{\rightarrow}}^{\text{Leaf}}\quad
\overbrace{\mathrm{M}}^{\text{Fork}}\quad
\overbrace{\binom{2}{\leftarrow}\mathrm{L}\binom{1}{\leftarrow}\mathrm{L}\ }^{\text{Trunk}}\quad
\overbrace{\binom{2}{\leftarrow}\mathrm{L}}^{\text{Leaf}}\)
}}
\put(0,-20){\makebox(0,0)[cc]{Transfer kernel types (TKT):}}
\put(0,-20.5){\makebox(0,0)[cc]{L: \(\varkappa=(12222222)\), M: \(\varkappa_0=(02222222)\)}}

\put(-1.2,-0.2){\framebox(0.4,0.4){}}
\put(-0.9,0.2){\makebox(0,0)[lb]{\(\langle 2\rangle\)}}
\put(-1,0){\line(0,-1){2}}
\put(-1,-2){\circle*{0.2}}
\put(-0.9,-1.8){\makebox(0,0)[lb]{\(\langle 3\rangle\)}}

\put(-1,-2){\vector(0,-1){2}}
\put(-1.2,-3.3){\makebox(0,0)[rc]{infinite}}
\put(-1.2,-3.8){\makebox(0,0)[rc]{mainline}}
\put(-1.2,-4.3){\makebox(0,0)[rc]{\(\mathcal{T}^1\langle 49,2\rangle\)}}

\put(-0.5,-2.8){\makebox(0,0)[lb]{metabelian bifurcation}}
\put(-1,-2){\line(1,-4){1}}

\put(0,-6){\line(0,-1){4}}
\multiput(0,-6)(0,-2){3}{\circle{0.2}}
\put(0.1,-5.8){\makebox(0,0)[lb]{\(\langle 4\rangle\)}}
\put(-0.1,-7.7){\makebox(0,0)[rb]{fork}}
\put(0.1,-7.8){\makebox(0,0)[lb]{\(\langle 708\rangle\)}}
\put(0.1,-9.8){\makebox(0,0)[lb]{\(\langle 2004\rangle\)}}
\put(0,-8){\line(-1,-2){1}}
\put(0,-8){\line(-1,-1){2}}
\multiput(-1,-10)(-1,0){2}{\circle{0.2}}
\put(-1,-10.2){\makebox(0,0)[ct]{\(\langle 1990\rangle\)}}
\put(-2,-10.2){\makebox(0,0)[ct]{\(\langle 1991\rangle\)}}

\put(-1.5,-10){\oval(2.5,1.2)}
\put(-3.8,-8.5){\makebox(0,0)[lc]{\(S/S^{\prime\prime}\simeq\)}}
\put(-4.8,-9){\makebox(0,0)[lc]{\(M=\mathrm{Gal}(\mathrm{F}_7^2(K)/K)\)}}

\put(0,-10){\vector(0,-1){2}}
\put(-0.2,-11.3){\makebox(0,0)[rc]{infinite}}
\put(-0.2,-11.8){\makebox(0,0)[rc]{mainline}}
\put(-0.2,-12.3){\makebox(0,0)[rc]{\(\mathcal{T}^2\langle 16807,4\rangle\)}}

\put(2.5,-8.8){\makebox(0,0)[lb]{first non-metabelian bifurcation}}
\put(0,-8){\line(1,-4){1}}
\put(0,-8){\line(4,-5){2.5}}
\put(0,-8){\line(3,-1){12}}

\multiput(0.9,-12.1)(0.5,0){7}{\framebox(0.2,0.2){}}
\put(0.5,-11.8){\makebox(0,0)[cb]{\(2;\)}}
\put(1,-11.8){\makebox(0,0)[cb]{\(113\)}}
\put(1.5,-11.8){\makebox(0,0)[cb]{\(96\)}}
\put(2,-11.8){\makebox(0,0)[cb]{\(79\)}}
\put(2.5,-11.8){\makebox(0,0)[cb]{\(62\)}}
\put(3,-11.8){\makebox(0,0)[cb]{\(45\)}}
\put(3.5,-11.8){\makebox(0,0)[cb]{\(28\)}}
\put(4,-11.8){\makebox(0,0)[cb]{\(11\)}}
\put(4.5,-12.8){\makebox(0,0)[cc]{\(\overbrace{\qquad\qquad\qquad 14\text{ times similarly } \qquad\qquad\qquad}\)}}
\multiput(4.9,-12.1)(0.5,0){7}{\framebox(0.2,0.2){}}
\put(5,-11.8){\makebox(0,0)[cb]{\(106\)}}
\put(5.5,-11.8){\makebox(0,0)[cb]{\(89\)}}
\put(6,-11.8){\makebox(0,0)[cb]{\(72\)}}
\put(6.5,-11.8){\makebox(0,0)[cb]{\(55\)}}
\put(7,-11.8){\makebox(0,0)[cb]{\(38\)}}
\put(7.5,-11.8){\makebox(0,0)[cb]{\(21\)}}
\put(8,-11.8){\makebox(0,0)[cb]{\(4\)}}
\multiput(8.95,-12.05)(0.5,0){7}{\framebox(0.1,0.1){}}
\put(9,-11.8){\makebox(0,0)[cb]{\(103\)}}
\put(9.5,-11.8){\makebox(0,0)[cb]{\(86\)}}
\put(10,-11.8){\makebox(0,0)[cb]{\(69\)}}
\put(10.5,-11.8){\makebox(0,0)[cb]{\(52\)}}
\put(11,-11.8){\makebox(0,0)[cb]{\(35\)}}
\put(11.5,-11.8){\makebox(0,0)[cb]{\(18\)}}
\put(12,-11.8){\makebox(0,0)[cb]{\(1\)}}

\multiput(1,-12)(4,0){2}{\line(0,-1){2}}
\multiput(1.5,-12)(4,0){2}{\line(0,-1){2}}
\multiput(2,-12)(4,0){2}{\line(0,-1){2}}
\multiput(2.5,-12)(4,0){2}{\line(0,-1){2}}
\multiput(3,-12)(4,0){2}{\line(0,-1){2}}
\multiput(3.5,-12)(4,0){2}{\line(0,-1){2}}
\multiput(4,-12)(4,0){2}{\line(0,-1){2}}
\multiput(2.4,-14.1)(4.5,0){1}{\framebox(0.2,0.2){}}
\multiput(2.5,-14.2)(4.5,0){1}{\makebox(0,0)[rt]{\(1;4\)}}

\put(4.5,-14.8){\makebox(0,0)[lb]{second non-metabelian bifurcation}}
\put(2.5,-14){\line(3,-4){3}}
\put(2.5,-14){\line(2,-1){7.8}}

\multiput(5.4,-18.1)(0.8,0){7}{\framebox(0.2,0.2){\vrule height 3pt width 3pt}}
\put(5.5,-18.2){\makebox(0,0)[ct]{\(2;1\)}}
\put(6.3,-18.2){\makebox(0,0)[ct]{\(2;2\)}}
\put(7.1,-18.2){\makebox(0,0)[ct]{\(2;3\)}}
\put(7.9,-18.2){\makebox(0,0)[ct]{\(2;4\)}}
\put(8.7,-18.2){\makebox(0,0)[ct]{\(2;5\)}}
\put(9.5,-18.2){\makebox(0,0)[ct]{\(2;6\)}}
\put(10.3,-18.2){\makebox(0,0)[ct]{\(2;7\)}}
\put(7.9,-19){\makebox(0,0)[cc]{\(7\) Schur \(\sigma\)-groups with \(\mathrm{dl}=3\)}}

\put(8,-18.5){\oval(7,1.5)}
\put(12,-19.8){\makebox(0,0)[rc]{\(S=\mathrm{Gal}(\mathrm{F}_7^\infty(K)/K)\)}}

}

\end{picture}

\end{figure}

\noindent
We apply the strategy of pattern recognition via Artin transfers
\cite{Ma2020}
to the situation described in Theorem
\ref{thm:7ClassTowers}.
The number theoretic Artin pattern
\(\mathrm{AP}(K)=((\ker(T_{K,E_i}))_{1\le i\le 8},(\mathrm{Cl}_7(E_i))_{1\le i\le 8})\)
of the class extension homomorphisms
\(T_{K,E_i}:\,\mathrm{Cl}_7(K)\to\mathrm{Cl}_7(E_i)\),
\(\mathfrak{a}\mathcal{P}_K\mapsto(\mathfrak{a}\mathcal{O}_{E_i})\mathcal{P}_{E_i}\),
from \(K\) to its eight unramified cyclic septic extension fields \(E_i\)
is interpreted as group theoretic Artin pattern
\(\mathrm{AP}(G)=((\ker(T_{G,H_i}))_{1\le i\le 8},(H_i/H_i^\prime)_{1\le i\le 8})\)
of the Artin transfer homomorphisms
\(T_{G,H_i}:\,G/G^\prime\to H_i/H_i^\prime\)
from any of the \(n\)-th \(7\)-class groups \(G=\mathrm{Gal}(\mathrm{F}_7^n(K)/K)\) of \(K\),
with \(n\ge 2\),
to its eight maximal subgroups \(H_i\) of index \(7\)
\cite{Ar1,Ar2,Ma2015c,Ma2016d,Ma2016e}.
The following proof is visualized by Figure
\ref{fig:SchurSigmaLogOrd11}.


We begin our search for the desired transfer kernel type (TKT) \(\varkappa\sim (1,2,2,2,2,2,2,2)\)
at the abelian root \(\mathrm{SmallGroup}(49,2)\simeq C_7\times C_7\)
of the tree of all finite \(7\)-groups \(G\) with abelianization \(G/G^\prime\) of type \((7,7)\).
Groups \(G\) with coclass \(\mathrm{cc}(G)=1\) in the coclass tree \(\mathcal{T}^1\langle 49,2\rangle\)
are discouraged, since they possess at least \(7\) total transfer kernels,
that is, \(\varkappa\sim (\ast,0,0,0,0,0,0,0)\) with a place holder \(\ast\),
but imaginary quadratic fields with \(\varrho_7=2\) cannot have total capitulation
\cite{Ma1991,Ma2016b,Ma2016c}.

Thus we follow the edges with step size \(s=2\) at the metabelian bifurcation of \(\mathrm{SmallGroup}(343,3)\).
Many of these immediate descendants have permutations as TKT,
but the search can be restricted to the tree of \(\mathrm{SmallGroup}(16\,807,4)\)
with TKT \(\varkappa\sim (0,2,2,2,2,2,2,2)\),
since this type is able to shrink to \((1,2,2,2,2,2,2,2)\) with two fixed points.
Note that the trees of \(\mathrm{SmallGroup}(16\,807,i)\) with \(5\le i\le 6\) are irrelevant,
since the TKT \(\varkappa\sim (0,1,1,1,1,1,1,1)\)
can only shrink to the type \((1,1,1,1,1,1,1,1)\) with a single fixed point
\cite{Bg,Ma2016a}.

The possible metabelian second \(7\)-class groups \(\mathrm{Gal}(\mathrm{F}_7^2(K)/K)\) of \(K\)
are easy to find as leaves of the second branch of the coclass tree \(\mathcal{T}^2\langle 16\,807,4\rangle\).
With respect to the branch root, they are given by the identifiers
\(\mathrm{SmallGroup}(823\,543,1990)=\mathrm{SmallGroup}(117\,649,708)-\#1;11\) and
\(\mathrm{SmallGroup}(823\,543,1991)=\mathrm{SmallGroup}(117\,649,708)-\#1;18\).
They are not unique with TKT \((1,2,2,2,2,2,2,2)\),
but they are unique in combination with the required transfer target type (TTT)
\(\tau\sim (1^5,1^3,21,21,21,21,21,21)\).
Since they have relation rank \(d_2=3>2=d_1\), their presentation is not balanced,
and the \(7\)-class field tower must have length \(\ell_7(K)\ge 3\)
\cite{Sh}.

In order to find the non-metabelian Schur \(\sigma\)-Galois group \(\mathrm{Gal}(\mathrm{F}_7^\infty(K)/K)\)
of the maximal unramified pro-\(7\) extension of \(K\),
we follow the paths with step size \(s=2\) at the non-metabelian bifurcation of \(\mathrm{SmallGroup}(117\,649,708)\)
with nuclear rank \(\nu=2\),
similarly as demonstrated for \(p=3\) in
\cite{BuMa}.
In contrast to \(p=3\), none of the immediate descendants of order \(7^8\) is a Schur \(\sigma\)-group.
The siblings \(G=\mathrm{SmallGroup}(117\,649,708)-\#2;\ell\) with

\begin{itemize}
\item
\(\ell\in\lbrace 1+17v\mid 0\le v\le 6\rbrace\) have type \((0,2,2,2,2,2,2,2)\),
\item
\(\ell\in\lbrace 4+17v\mid 0\le v\le 6\rbrace\) have type \((1,2,2,2,2,2,2,2)\) and \(G/G^{\prime\prime}\simeq\langle 823\,543,1990\rangle\),
\item
\(\ell\in\lbrace 11+17v\mid 0\le v\le 6\rbrace\) have type \((1,2,2,2,2,2,2,2)\) and \(G/G^{\prime\prime}\simeq\langle 823\,543,1991\rangle\).
\end{itemize}

The danger that descendants of the groups with \(\ell\in\lbrace 1+17v\mid 0\le v\le 6\rbrace\)
give rise to further relevant Schur \(\sigma\)-groups is eliminated by the monotony principle
\cite{Ma2016a},
since the polarization of the TTT \(\tau\) increases
from \(\tau_1=1^5\) with rank \(5\) to \(\tau_1=1^6\) with rank \(6\).

After a second non-metabelian bifurcation at the unique \(\sigma\)-descendant
of a group with \(\ell\in\lbrace 4+17v\mid 0\le v\le 6\rbrace\)
or \(\ell\in\lbrace 11+17v\mid 0\le v\le 6\rbrace\),
we arrive at the desired Schur \(\sigma\)-groups \(S\) with \(\#S=7^{11}\) and \(\mathrm{dl}(S)=3\).
These are \(49\) candidates for each of the two metabelianizations.
Their identifiers are given by \(\mathrm{SmallGroup}(117\,649,708)-\#2;\ell-\#1,k-\#2;j\)
as in Formula
\eqref{eqn:7ClassTowers}.


In Table
\ref{tbl:RootOrd16807Id4},
invariants for the pre-period and period
of the tree with root \(\langle 7^5,4\rangle\)
and its descendants of fixed coclass \(\mathrm{cc}=2\)
on the mainline or with type \(\varkappa\) in Theorem
\ref{thm:7ClassTowers}
are given.
The identifier \(\mathrm{id}\) of each vertex
is either the absolute number in the SmallGroups database
\cite{BEO1,BEO2}
or the relative counter with respect to the parent,
in the notation of the ANUPQ package
\cite{GNO}.
The order \(\mathrm{ord}\) and nilpotency class \(\mathrm{cl}\)
are arranged in ascending order.
The Artin pattern \(\mathrm{AP}=(\varkappa,\tau)\)
consists of the \textit{stabilization}
\cite{Ma2016a}

\begin{equation}
\label{eqn:Stabilization7}
\varkappa\sim(\ast,2,2,2,2,2,2,2), \quad \tau\sim(\ast,1^3,21,21,21,21,21,21)
\end{equation}

\noindent
with a placeholder \(\ast\) for the \textit{polarization} \((\varkappa_1,\tau_1)\).
The stabilization remains constant along the entire coclass tree,
whereas the polarization varies from branch to branch
with strictly increasing \(\tau_1\).
Therefore, it suffices to list the polarization.
The graph theoretic structure of the tree is determined by
the \textit{nuclear rank} \(\nu\) and the descendant numbers \(N_1/C_1\) with step size \(s=1\),
respectively \(N_2/C_2\) with step size \(s=2\), when \(\nu=2\).
The latter phenomenon, which is called \textit{bifurcation} to bigger coclass,
occurs three times in the irregular pre-period of length \(6\).
Starting with the first primitive period, also of length \(6\),
the value of the nuclear rank \(\nu=1\) remains settled,
and no further information is required in order to describe the entire infinite tree.
The \(p\)-multiplicator rank \(\mu\) coincides with the \textit{relation rank} \(d_2\)
in the Shafarevich theorem
\cite{Sh}.


\renewcommand{\arraystretch}{1.1}

\begin{table}[ht]
\caption{Structure and periodicity of the coclass tree \(\mathcal{T}^2\langle 16\,807,4\rangle\)}
\label{tbl:RootOrd16807Id4}
\begin{center}
\begin{tabular}{|rl|rl|r|cc|cc|lc|}
\hline
\multicolumn{2}{|c|}{id} & \multicolumn{2}{c|}{ord}                 & cl     & \(\varkappa_1\) & \(\tau_1\) & \(\nu\) & \(\mu\) & \(N_1/C_1\) & \(N_2/C_2\) \\
\hline
     \(4\) &             &                \(16\,807\) &    \(=7^5\) &  \(3\) &           \(0\) &    \(1^3\) &   \(1\) &   \(3\) & \(11/11\)   &             \\
   \(708\) &  \(-\#1;1\) &               \(117\,649\) &    \(=7^6\) &  \(4\) &           \(0\) &    \(1^4\) &   \(2\) &   \(4\) & \(38/3\)    & \(119/119\) \\
  \(2004\) &  \(-\#1;8\) &               \(823\,543\) &    \(=7^7\) &  \(5\) &           \(0\) &    \(1^5\) &   \(1\) &   \(4\) & \(32/14\)   &             \\
  \(1990\) & \(-\#1;11\) &               \(823\,543\) &    \(=7^7\) &  \(5\) &           \(1\) &    \(1^5\) &   \(0\) &   \(3\) &             &             \\ 
  \(1991\) & \(-\#1;18\) &               \(823\,543\) &    \(=7^7\) &  \(5\) &           \(1\) &    \(1^5\) &   \(0\) &   \(3\) &             &             \\ 
           &  \(-\#1;1\) &            \(5\,764\,801\) &    \(=7^8\) &  \(6\) &           \(0\) &    \(1^6\) &   \(2\) &   \(5\) & \(52/9\)    & \(805/217\) \\
           & \(-\#1;10\) &           \(40\,353\,607\) &    \(=7^9\) &  \(7\) &           \(0\) &   \(21^5\) &   \(1\) &   \(5\) & \(25/13\)   &             \\
           &  \(-\#1;2\) &          \(282\,475\,249\) & \(=7^{10}\) &  \(8\) &           \(0\) & \(2^21^4\) &   \(2\) &   \(6\) & \(80/11\)   & \(833/119\) \\
\hline
           &  \(-\#1;4\) &       \(1\,977\,326\,743\) & \(=7^{11}\) &  \(9\) &           \(0\) & \(2^31^3\) &   \(1\) &   \(6\) & \(28/6\)    &             \\
           &  \(-\#1;3\) &      \(13\,841\,287\,201\) & \(=7^{12}\) & \(10\) &           \(0\) & \(2^41^2\) &   \(1\) &   \(6\) & \(51/11\)   &             \\
           &  \(-\#1;2\) &      \(96\,889\,010\,407\) & \(=7^{13}\) & \(11\) &           \(0\) &   \(2^51\) &   \(1\) &   \(6\) & \(40/6\)    &             \\
           &  \(-\#1;2\) &     \(678\,223\,072\,849\) & \(=7^{14}\) & \(12\) &           \(0\) &    \(2^6\) &   \(1\) &   \(6\) & \(55/11\)   &             \\
           &  \(-\#1;3\) &  \(4\,747\,561\,509\,943\) & \(=7^{15}\) & \(13\) &           \(0\) &   \(32^5\) &   \(1\) &   \(6\) & \(26/6\)    &             \\
           &  \(-\#1;2\) & \(33\,232\,930\,569\,601\) & \(=7^{16}\) & \(14\) &           \(0\) & \(3^22^4\) &   \(1\) &   \(6\) & \(79/11\)   &             \\
\hline
\end{tabular}
\end{center}
\end{table}


In view of Table
\ref{tbl:RootOrd16807Id4}
it is obvious that 
\(\mathrm{SmallGroup}(823\,543,1990)=\mathrm{SmallGroup}(117\,649,708)-\#1;11\)
and
\(\mathrm{SmallGroup}(823\,543,1991)=\mathrm{SmallGroup}(117\,649,708)-\#1;18\)
are the only two possibilities for the metabelianization \(M\simeq S/S^{\prime\prime}\)
of the Galois group \(S\) of the \(7\)-class field tower of any algebraic number field \(K\)
satisfying the conditions in Theorem
\ref{thm:7ClassTowers},
because firstly the coclass tree \(\mathcal{T}^2\langle 7^5,4\rangle\)
with metabelian main line is unique with stabilization in Formula
\eqref{eqn:Stabilization7},
and secondly,
although there occur vertices with \(\varkappa\sim (1,2,2,2,2,2,2,2)\) on every branch of the tree,
the mentioned groups on branch \(\mathcal{B}(2)\) with root \(\mathrm{SmallGroup}(117\,649,708)\)
are the only two which possess the required polarization \(\tau_1=1^5\).
Since their relation rank is given by \(d_2=\mu=3\), bigger than the required \(d_2=2\),
a first consequence is the lower bound \(\ell_7(K)\ge 3\) for the length of the \(7\)-class tower of \(K\).


In Table
\ref{tbl:RootOrd117649Id708},
we give an exemplary \textit{root path}
to one of the non-metabelian Schur \(\sigma\)-groups \(S\simeq\mathrm{Gal}(\mathrm{F}_7^\infty(K)/K)\)
and corresponding invariants.
Figure
\ref{fig:SchurSigmaLogOrd11}
illustrates Table
\ref{tbl:RootOrd16807Id4}
and Table
\ref{tbl:RootOrd117649Id708}.


\renewcommand{\arraystretch}{1.1}

\begin{table}[ht]
\caption{Exemplary root path to \(S\simeq\mathrm{Gal}(\mathrm{F}_7^\infty(K)/K)\) with \(\mathrm{lo}=11\)}
\label{tbl:RootOrd117649Id708}
\begin{center}
\begin{tabular}{|rl|rl|r|r|cc|cc|lc|}
\hline
\multicolumn{2}{|c|}{id} & \multicolumn{2}{c|}{ord}                 & cl    & dl    & \(\varkappa_1\) & \(\tau_1\) & \(\nu\) & \(\mu\) & \(N_1/C_1\) & \(N_2/C_2\) \\
\hline
   \(708\) &             &               \(117\,649\) &    \(=7^6\) & \(4\) & \(2\) &           \(0\) &    \(1^4\) &   \(2\) &   \(4\) & \(38/3\)    & \(119/119\) \\
           & \(-\#2;62\) &            \(5\,764\,801\) &    \(=7^8\) & \(5\) & \(3\) &           \(1\) &    \(1^5\) &   \(1\) &   \(3\) &  \(4/4\)    &             \\
           &  \(-\#1;4\) &           \(40\,353\,607\) &    \(=7^9\) & \(6\) & \(3\) &           \(1\) &    \(1^5\) &   \(2\) &   \(4\) & \(14/0\)    &   \(7/0\)   \\
           &  \(-\#2;1\) &       \(1\,977\,326\,743\) & \(=7^{11}\) & \(7\) & \(3\) &           \(1\) &    \(1^5\) &   \(0\) &   \(2\) &             &             \\
\hline
\end{tabular}
\end{center}
\end{table}



\section{\(5\)-class field towers with length at least three}
\label{s:5ClassTowers}

\noindent
The analogue of the results in \S\
\ref{s:7ClassTowers}
for \(p=5\) is astonishing and disappointing for two reasons:

\begin{itemize}
\item
Firstly, the smallest Schur \(\sigma\)-groups with desired Artin pattern
set in at the annoying logarithmic order \(14\) instead of \(11\),
\item
and secondly, even more annoying,
there is not a finite number of such Schur \(\sigma\)-groups
but rather an infinitude of them with unbounded derived length,
similarly as known from the situations with sporadic metabelianizations
\cite{BaBu}.
\end{itemize}


\begin{theorem}
\label{thm:5ClassTowers} 
(Sufficient criterion for a \(5\)-class field tower with at least three stages) \\
An imaginary quadratic number field \(K=\mathbb{Q}(\sqrt{d})\)
with fundamental discriminant \(d<0\),
elementary bicyclic \(5\)-class group \(\mathrm{Cl}_5(K)\simeq C_5\times C_5\),
and Artin pattern \(\mathrm{AP}(K)=(\varkappa(K),\tau(K))\)
with capitulation type \(\varkappa(K)\sim (1,2,2,2,2,2)\)
in the six unramified cyclic septic extensions \((E_i)_{1\le i\le 6}\) of \(K\)
and logarithmic abelian type invariants \(\tau(K)\sim (21^3,1^3,21,21,21,21)\)
of the \(5\)-class groups \((\mathrm{Cl}_5(E_i))_{1\le i\le 6}\)
possesses a finite \(5\)-class field tower
\(K<\mathrm{F}_5^1(K)<\mathrm{F}_5^2(K)<\mathrm{F}_5^3(K)\le\mathrm{F}_5^\infty(K)\)
with length \(\ell_5(K)\ge 3\).
\end{theorem}


\noindent
Now we supplement the purely arithmetical information of Theorem
\ref{thm:5ClassTowers}
with group theoretical data.
We use the notation of the SmallGroups database
\cite{BEO1,BEO2}
and the ANUPQ package
\cite{GNO}.


\begin{corollary}
\label{cor:5ClassTowers} 
Under the conditions of Theorem
\ref{thm:5ClassTowers},
there are only two possibilities for the metabelian second \(5\)-class group
\(M=\mathrm{Gal}(\mathrm{F}_5^2(K)/K)\) of \(K\), namely
\(M\simeq\mathrm{SmallGroup}(78\,125,m)\) with either \(m=889\) or \(m=890\).
Both share the invariants \(\mathrm{lo}=7\), \(\mathrm{cl}=5\), \(\mathrm{cc}=2\), \(d_2=3\).
See Figure
\ref{fig:SchurSigmaLogOrd14}.

There are \textbf{infinitely many} candidates for the non-metabelian Schur \(\sigma\)-Galois group 
\(S=\mathrm{Gal}(\mathrm{F}_5^\infty(K)/K)\)
of the \(5\)-class field tower of \(K\).
In dependence on the metabelianization \(M\simeq S/S^{\prime\prime}\),
the \(200\) candidates with smallest possible order \(5^{14}\) and derived length \(3\) are given by

\begin{equation}
\label{eqn:5ClassTowers}
\begin{aligned}
S &\simeq\mathrm{SmallGroup}(15\,625,564)-\#2;\ell-\#1,k-\#2;j-\#1,i-\#2;h \\
\text{ with } &1\le h\le 5 \text{ and} \\
(\ell,k) &=(21,1),\ (j,i)\in\lbrace (1,11),(2,3),(3,12),(4,13),(5,8)\rbrace \text{ or} \\
(\ell,k) &=(22,3),\ (j,i)\in\lbrace (1,6),(2,13),(3,2),(4,13),(5,13)\rbrace \text{ or} \\
(\ell,k) &=(23,3),\ (j,i)\in\lbrace (1,13),(2,6),(3,13),(4,2),(5,13)\rbrace \text{ or} \\
(\ell,k) &=(25,2),\ (j,i)\in\lbrace (1,3),(2,11),(3,13),(4,7),(5,13)\rbrace, \text{ if } m=889, \\
(\ell,k) &=(26,1),\ (j,i)\in\lbrace (1,13),(2,12),(3,13),(4,3),(5,6)\rbrace \text{ or} \\
(\ell,k) &=(27,3),\ (j,i)\in\lbrace (1,7),(2,13),(3,13),(4,11),(5,3)\rbrace \text{ or} \\
(\ell,k) &=(29,3),\ (j,i)\in\lbrace (1,11),(2,3),(3,7),(4,13),(5,13)\rbrace \text{ or} \\
(\ell,k) &=(30,2),\ (j,i)\in\lbrace (1,12),(2,13),(3,3),(4,6),(5,13)\rbrace, \text{ if } m=890. \\
\end{aligned}
\end{equation}

\noindent
The latter share the invariants \(\mathrm{lo}=14\), \(\mathrm{cl}=9\), \(\mathrm{cc}=5\),
and \(d_2=2=d_1=\dim_{\mathbb{F}_5}(\mathrm{H}^1(S,\mathbb{F}_5))=\varrho_5\).

However, there also exist \(25\,000\) candidates with order \(5^{21}\) and derived length \(4\)
among the descendants of the two groups \(\mathrm{SmallGroup}(15\,625,564)-\#2;\ell\) with \(\ell\in\lbrace 24,28\rbrace\).
See Figure
\ref{fig:SchurSigmaLogOrd21}
\end{corollary}


\begin{proof}
The next section is devoted to the rigorous justification of Theorem
\ref{thm:5ClassTowers}
and Corollary
\ref{cor:5ClassTowers}.
\end{proof}


\begin{example}
\label{exm:5ClassTowers}
The absolutely smallest discriminant of an imaginary quadratic number field \(K=\mathbb{Q}(\sqrt{d})\)
which satisfies the conditions in Theorem
\ref{thm:5ClassTowers}
is \(d=-90\,868\).
The \(5\)-class field tower over this field has at least three stages,
but the precise length \(\ell_5(K)\ge 3\) remains unknown.
\end{example} 



\section{\(5\)-groups connected with \(5\)-class towers}
\label{s:5Groups}

\noindent
We employ the strategy of pattern recognition via Artin transfers
\cite{Ma2020}
to the situation described in Theorem
\ref{thm:5ClassTowers}.
The number theoretic Artin pattern
\(\mathrm{AP}(K)=((\ker(T_{K,E_i}))_{1\le i\le 6},(\mathrm{Cl}_5(E_i))_{1\le i\le 6})\)
of the class extension homomorphisms
\(T_{K,E_i}:\,\mathrm{Cl}_5(K)\to\mathrm{Cl}_5(E_i)\),
\(\mathfrak{a}\mathcal{P}_K\mapsto(\mathfrak{a}\mathcal{O}_{E_i})\mathcal{P}_{E_i}\),
from \(K\) to its six unramified cyclic quintic extension fields \(E_i\)
is interpreted as group theoretic Artin pattern
\(\mathrm{AP}(G)=((\ker(T_{G,H_i}))_{1\le i\le 6},(H_i/H_i^\prime)_{1\le i\le 6})\)
of the Artin transfer homomorphisms
\(T_{G,H_i}:\,G/G^\prime\to H_i/H_i^\prime\)
from any of the \(n\)-th \(5\)-class groups \(G=\mathrm{Gal}(\mathrm{F}_5^n(K)/K)\) of \(K\),
with \(n\ge 2\),
to its six maximal subgroups \(H_i\) of index \(5\)
\cite{Ar1,Ar2,Ma2015c,Ma2016d,Ma2016e}.
The following proof is visualized by Fig.
\ref{fig:SchurSigmaLogOrd14}.


We begin our search for the desired transfer kernel type (TKT) \(\varkappa\sim (1,2,2,2,2,2)\)
at the abelian root \(\mathrm{SmallGroup}(25,2)\simeq C_5\times C_5\)
of the tree of all finite \(5\)-groups \(G\) with abelianization \(G/G^\prime\) of type \((5,5)\).
Groups \(G\) with coclass \(\mathrm{cc}(G)=1\) in the coclass tree \(\mathcal{T}^1\langle 25,2\rangle\)
are discouraged, since they possess at least \(5\) total transfer kernels,
that is, \(\varkappa\sim (\ast,0,0,0,0,0)\) with a place holder \(\ast\),
but imaginary quadratic fields with \(\varrho_5=2\) cannot have total capitulation
\cite{Ma1991,Ma2016b,Ma2016c}.

Thus we follow the edges with step size \(s=2\) at the metabelian bifurcation of \(\mathrm{SmallGroup}(125,3)\).
Many of these immediate descendants have permutations as TKT,
but the search can be restricted to the tree of \(\mathrm{SmallGroup}(3\,125,4)\)
with TKT \(\varkappa\sim (0,2,2,2,2,2)\),
since this type is able to shrink to \((1,2,2,2,2,2)\) with two fixed points.
Note that the trees of \(\mathrm{SmallGroup}(3\,125,i)\) with \(5\le i\le 6\) are irrelevant,
since the TKT \(\varkappa\sim (0,1,1,1,1,1)\)
can only shrink to the type \((1,1,1,1,1,1)\) with a single fixed point
\cite{Bg,Ma2016a}.

The possible metabelian second \(5\)-class groups \(\mathrm{Gal}(\mathrm{F}_5^2(K)/K)\) of \(K\)
are easy to find as leaves of the second branch of the coclass tree \(\mathcal{T}^2\langle 3\,125,4\rangle\).
With respect to the branch root, they are given by the identifiers
\(\langle 78\,125,889\rangle=\mathrm{SmallGroup}(15\,625,564)-\#1;5\) and
\(\langle 78\,125,890\rangle=\mathrm{SmallGroup}(15\,625,564)-\#1;6\).
They are not unique with TKT \((1,2,2,2,2,2)\),
but they are unique in combination with the required transfer target type (TTT)
\(\tau\sim (21^3,1^3,21,21,21,21)\).
Since they have relation rank \(d_2=3>2=d_1\), their presentation is not balanced,
and the \(5\)-class field tower must have length \(\ell_5(K)\ge 3\)
\cite{Sh}.

In order to find the non-metabelian Schur \(\sigma\)-Galois group \(\mathrm{Gal}(\mathrm{F}_5^\infty(K)/K)\)
of the maximal unramified pro-\(5\) extension of \(K\),
we follow the paths with step size \(s=2\) at the non-metabelian bifurcation of \(\mathrm{SmallGroup}(15\,625,564)\)
with nuclear rank \(\nu=2\),
similarly as demonstrated for \(p=3\) in
\cite{BuMa}.
In contrast to \(p=3\), none of the immediate descendants of order \(5^8\) is a Schur \(\sigma\)-group.
The siblings \(G=\mathrm{SmallGroup}(15\,625,564)-\#2;\ell\) with

\begin{itemize}
\item
\(16\le\ell\le 20\) have type \((0,2,2,2,2,2)\),
\item
\(21\le\ell\le 25\) have type \((1,2,2,2,2,2)\) and \(G/G^{\prime\prime}\simeq\langle 78\,125,889\rangle\),
\item
\(26\le\ell\le 30\) have type \((1,2,2,2,2,2)\) and \(G/G^{\prime\prime}\simeq\langle 78\,125,890\rangle\).
\end{itemize}

The danger that descendants of the groups with \(16\le\ell\le 20\)
give rise to further relevant Schur \(\sigma\)-groups is eliminated by the monotony principle
\cite{Ma2016a},
since the polarization of the TTT \(\tau\) increases
from \(\tau_1=21^3\) to \(\tau_1=2^21^2\), both with rank \(4\).

After a second non-metabelian bifurcation at the unique \(\sigma\)-descendant
of a group with \(\ell\in\lbrace 21,22,23,25\rbrace\),
resp. \(\ell\in\lbrace 26,27,29,30\rbrace\),
and a third non-metabelian bifurcation at the unique \(\sigma\)-descendant
of \(\mathrm{SmallGroup}(15\,625,564)-\#2;\ell-\#1;k-\#2;j\) with \(1\le j\le 5\)
we arrive at the desired Schur \(\sigma\)-groups \(S\) with \(\#S=5^{14}\) and \(\mathrm{dl}(S)=3\).
These are \(100\) candidates for each of the two metabelianizations.
Their identifiers are given by \(\mathrm{SmallGroup}(15\,625,564)-\#2;\ell-\#1,k-\#2;j-\#1,i-\#2;h\)
as in Formula
\eqref{eqn:5ClassTowers}.

\begin{figure}[ht]
\caption{Exceptional \(5\)-class field tower groups of order \(5^{21}\) on the tree \(\mathcal{T}\langle 3\,125,4\rangle\)}
\label{fig:SchurSigmaLogOrd21}


\setlength{\unitlength}{0.7cm}
\begin{picture}(18,26.5)(-6,-25.5)

{\tiny

\put(-5,0.5){\makebox(0,0)[cb]{Order}}
\put(-5,0){\line(0,-1){26}}
\multiput(-5.1,0)(0,-2){14}{\line(1,0){0.2}}
\put(-5.2,0){\makebox(0,0)[rc]{\(390\,625\)}}
\put(-4.8,0){\makebox(0,0)[lc]{\(5^8\)}}
\put(-5.2,-2){\makebox(0,0)[rc]{\(1\,953\,125\)}}
\put(-4.8,-2){\makebox(0,0)[lc]{\(5^9\)}}
\put(-5.2,-4){\makebox(0,0)[rc]{\(9\,765\,625\)}}
\put(-4.8,-4){\makebox(0,0)[lc]{\(5^{10}\)}}
\put(-5.2,-6){\makebox(0,0)[rc]{\(48\,828\,125\)}}
\put(-4.8,-6){\makebox(0,0)[lc]{\(5^{11}\)}}
\put(-5.2,-8){\makebox(0,0)[rc]{\(244\,140\,625\)}}
\put(-4.8,-8){\makebox(0,0)[lc]{\(5^{12}\)}}
\put(-5.2,-10){\makebox(0,0)[rc]{\(1\,220\,703\,125\)}}
\put(-4.8,-10){\makebox(0,0)[lc]{\(5^{13}\)}}
\put(-5.2,-12){\makebox(0,0)[rc]{\(6\,103\,515\,625\)}}
\put(-4.8,-12){\makebox(0,0)[lc]{\(5^{14}\)}}
\put(-5.2,-14){\makebox(0,0)[rc]{\(30\,517\,578\,125\)}}
\put(-4.8,-14){\makebox(0,0)[lc]{\(5^{15}\)}}
\put(-5.2,-16){\makebox(0,0)[rc]{\(152\,587\,890\,625\)}}
\put(-4.8,-16){\makebox(0,0)[lc]{\(5^{16}\)}}
\put(-5.2,-18){\makebox(0,0)[rc]{\(762\,939\,453\,125\)}}
\put(-4.8,-18){\makebox(0,0)[lc]{\(5^{17}\)}}
\put(-5.2,-20){\makebox(0,0)[rc]{\(3\,814\,697\,265\,625\)}}
\put(-4.8,-20){\makebox(0,0)[lc]{\(5^{18}\)}}
\put(-5.2,-22){\makebox(0,0)[rc]{\(19\,073\,486\,328\,125\)}}
\put(-4.8,-22){\makebox(0,0)[lc]{\(5^{19}\)}}
\put(-5.2,-24){\makebox(0,0)[rc]{\(95\,367\,431\,640\,625\)}}
\put(-4.8,-24){\makebox(0,0)[lc]{\(5^{20}\)}}
\put(-5.2,-26){\makebox(0,0)[rc]{\(476\,837\,158\,203\,125\)}}
\put(-4.8,-26){\makebox(0,0)[lc]{\(5^{21}\)}}
\put(-5,-26){\vector(0,-1){2}}


\put(-3,-15){\makebox(0,0)[lc]{Legend:}}
\put(-3.1,-16.1){\framebox(0.2,0.2){}}
\put(-3.1,-16.6){\framebox(0.2,0.2){}}
\put(-3,-16.5){\makebox(0,0)[cc]{\(\times\)}}
\put(-3.1,-17.1){\framebox(0.2,0.2){\vrule height 3pt width 3pt}}
\put(-3.05,-17.55){\framebox(0.1,0.1){}}

\put(-2.5,-16){\makebox(0,0)[lc]{non-metabelian with TKT L}}
\put(-2.5,-16.5){\makebox(0,0)[lc]{exceptional non-metabelian with TKT L}}
\put(-2.5,-17){\makebox(0,0)[lc]{non-metabelian Schur \(\sigma\)-group with TKT L}}
\put(-2.5,-17.5){\makebox(0,0)[lc]{non-metabelian with TKT M}}

\multiput(-2.1,-0.1)(0.5,0){5}{\framebox(0.2,0.2){}}
\put(-2.5,0.2){\makebox(0,0)[cb]{\(2;\)}}
\put(-2,0.2){\makebox(0,0)[cb]{\(30\)}}
\put(-1.5,0.2){\makebox(0,0)[cb]{\(29\)}}
\put(-1,0.2){\makebox(0,0)[cb]{\(28\)}}
\put(-1,0){\makebox(0,0)[cc]{\(\times\)}}
\put(-0.5,0.2){\makebox(0,0)[cb]{\(27\)}}
\put(0,0.2){\makebox(0,0)[cb]{\(26\)}}
\put(1,-0.8){\makebox(0,0)[cc]{\(\overbrace{\qquad 2\text{ times similarly } \qquad}\)}}
\multiput(2.4,-0.1)(0.5,0){5}{\framebox(0.2,0.2){}}
\put(2.5,0.2){\makebox(0,0)[cb]{\(25\)}}
\put(3,0.2){\makebox(0,0)[cb]{\(24\)}}
\put(3,0){\makebox(0,0)[cc]{\(\times\)}}
\put(3.5,0.2){\makebox(0,0)[cb]{\(23\)}}
\put(4,0.2){\makebox(0,0)[cb]{\(22\)}}
\put(4.5,0.2){\makebox(0,0)[cb]{\(21\)}}
\multiput(6.95,-0.05)(0.5,0){5}{\framebox(0.1,0.1){}}
\put(7,0.2){\makebox(0,0)[cb]{\(20\)}}
\put(7.5,0.2){\makebox(0,0)[cb]{\(19\)}}
\put(8,0.2){\makebox(0,0)[cb]{\(18\)}}
\put(8.5,0.2){\makebox(0,0)[cb]{\(17\)}}
\put(9,0.2){\makebox(0,0)[cb]{\(16\)}}

\multiput(-1,0)(4.5,0){1}{\line(0,-1){2}}
\multiput(3,0)(4.5,0){1}{\line(0,-1){2}}
\multiput(2.9,-2.1)(4.5,0){1}{\framebox(0.2,0.2){}}
\multiput(2.9,-1.8)(4.5,0){1}{\makebox(0,0)[rb]{\(1;3\)}}

\put(5,-2.8){\makebox(0,0)[lb]{bifurcation}}
\put(3,-2){\line(1,-1){4}}
\put(3,-2){\line(2,-1){8}}
\multiput(6.9,-6.1)(4,0){2}{\framebox(0.2,0.2){}}
\put(7,-5.8){\makebox(0,0)[lb]{\(2;1\)}}
\put(7,-8.1){\makebox(0,0)[cc]{\(\overbrace{\qquad 20\text{ times similarly } \qquad}\)}}
\put(9,-5.8){\makebox(0,0)[cb]{\(\ldots\)}}
\put(11,-5.8){\makebox(0,0)[lb]{\(2;20\)}}

\put(6,-6.8){\makebox(0,0)[rb]{bifurcation}}
\multiput(7,-6)(4,0){2}{\line(-1,-1){4}}
\multiput(2.9,-10.1)(4.5,0){1}{\framebox(0.2,0.2){}}
\multiput(2.9,-9.8)(4.5,0){1}{\makebox(0,0)[rb]{\(2;306\)}}

\put(5,-10.8){\makebox(0,0)[lb]{trifurcation}}
\put(3,-10){\line(2,-3){4}}
\put(3,-10){\line(4,-3){8}}

\multiput(6.9,-16.1)(4,0){2}{\framebox(0.2,0.2){}}
\put(7,-15.8){\makebox(0,0)[lb]{\(3;1\)}}
\put(9,-16.8){\makebox(0,0)[cc]{\(\overbrace{\qquad 25\text{ times similarly } \qquad}\)}}
\put(9,-15.8){\makebox(0,0)[cb]{\(\ldots\)}}
\put(11,-15.8){\makebox(0,0)[lb]{\(3;25\)}}

\multiput(7,-16)(4,0){2}{\line(0,-1){2}}
\multiput(6.9,-18.1)(4.5,0){1}{\framebox(0.2,0.2){}}
\multiput(6.9,-17.8)(4.5,0){1}{\makebox(0,0)[rb]{\(1;11\)}}

\put(5,-18.8){\makebox(0,0)[rb]{bifurcation}}
\put(7,-18){\line(-1,-1){4}}
\put(7,-18){\line(-2,-1){8}}
\multiput(-1.1,-22.1)(1,0){5}{\framebox(0.2,0.2){}}
\put(-1,-21.8){\makebox(0,0)[cb]{\(2;5\)}}
\put(0,-21.8){\makebox(0,0)[cb]{\(2;4\)}}
\put(1,-22.8){\makebox(0,0)[cc]{\(\overbrace{\qquad 5\text{ times similarly } \qquad}\)}}
\put(1,-21.8){\makebox(0,0)[cb]{\(2;3\)}}
\put(2,-21.8){\makebox(0,0)[cb]{\(2;2\)}}
\put(3,-21.8){\makebox(0,0)[cb]{\(2;1\)}}

\multiput(3,-22)(-1,0){5}{\line(0,-1){2}}
\multiput(2.9,-24.1)(4.5,0){1}{\framebox(0.2,0.2){}}
\multiput(3,-23.8)(4.5,0){1}{\makebox(0,0)[rb]{\(1;3\)}}
\put(3,-24){\line(-1,-1){2}}
\put(3,-24){\line(1,-1){2}}

\multiput(0.9,-26.1)(1,0){5}{\framebox(0.2,0.2){\vrule height 3pt width 3pt}}
\put(1,-26.2){\makebox(0,0)[ct]{\(1;1\)}}
\put(2,-26.2){\makebox(0,0)[ct]{\(1;2\)}}
\put(3,-26.2){\makebox(0,0)[ct]{\(1;3\)}}
\put(4,-26.2){\makebox(0,0)[ct]{\(1;4\)}}
\put(5,-26.2){\makebox(0,0)[ct]{\(1;5\)}}
\put(3,-27){\makebox(0,0)[cc]{\(5\) Schur \(\sigma\)-groups with \(\mathrm{dl}=4\)}}

\put(3,-26.4){\oval(5.5,1.5)}
\put(0.3,-25.5){\makebox(0,0)[rc]{\(S=\mathrm{Gal}(\mathrm{F}_5^\infty(K)/K)\)}}

}

\end{picture}

\end{figure}
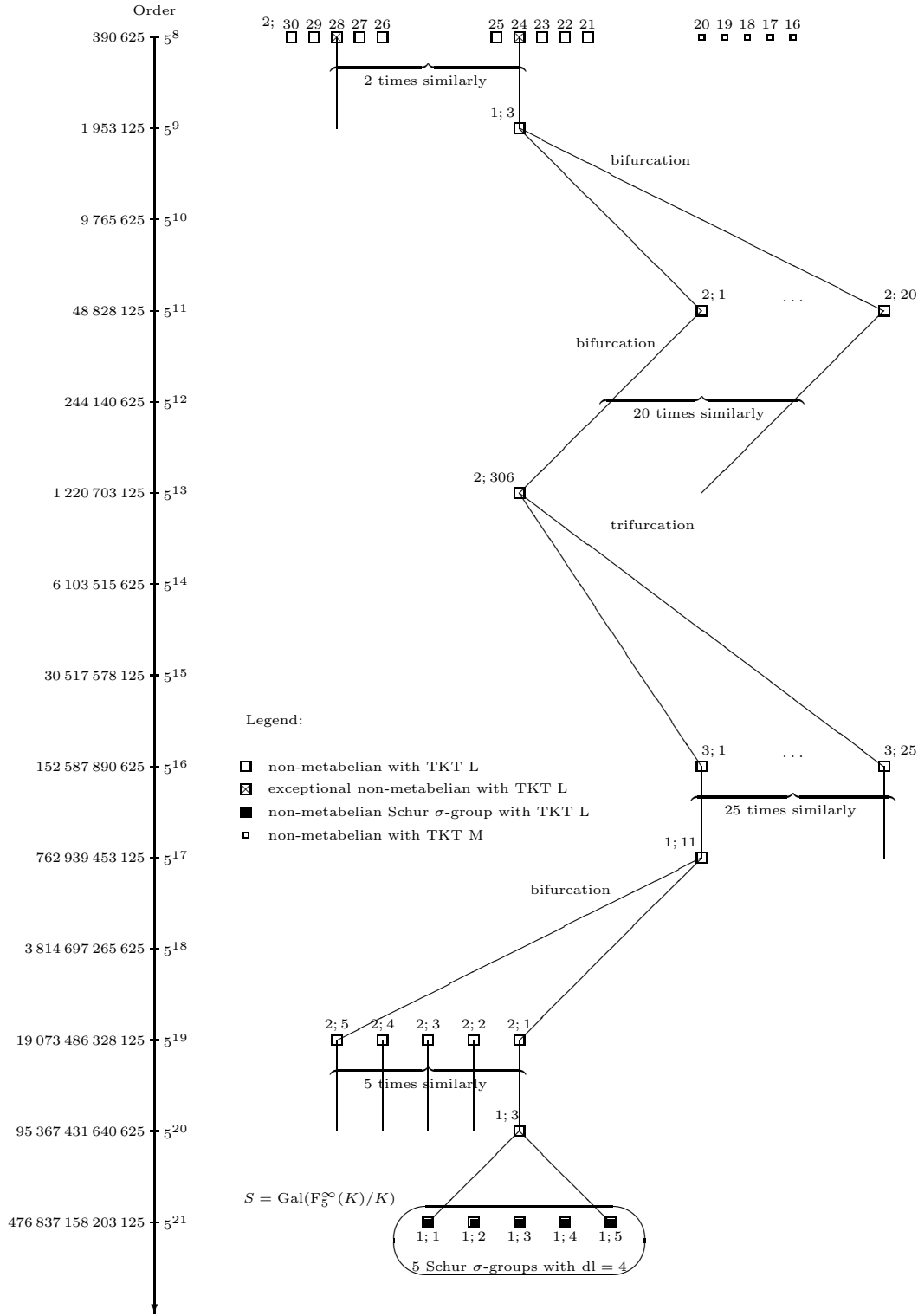

Unfortunately, the two groups with \(\ell\in\lbrace 24,28\rbrace\) are exceptional.
They destroy the warranty for a \(5\)-class field tower with precise length \(\ell_5(K)=3\).
See Table
\ref{tbl:RootOrd15625Id564Lo21}
and Figure
\ref{fig:SchurSigmaLogOrd21}.


In Table
\ref{tbl:RootOrd3125Id4},
invariants for the pre-period and period
of the tree with root \(\langle 5^5,4\rangle\)
and its descendants of fixed coclass \(\mathrm{cc}=2\)
on the mainline or with type \(\varkappa\) in Theorem
\ref{thm:5ClassTowers}
are given.
The identifier \(\mathrm{id}\) of each vertex
is either the absolute number in the SmallGroups database
\cite{BEO1,BEO2}
or the relative counter with respect to the parent,
in the notation of the ANUPQ package
\cite{GNO}.
The order \(\mathrm{ord}\) and nilpotency class \(\mathrm{cl}\)
are arranged in ascending order.
The Artin pattern \(\mathrm{AP}=(\varkappa,\tau)\)
consists of the \textit{stabilization}
\cite{Ma2016a}

\begin{equation}
\label{eqn:Stabilization5}
\varkappa\sim(\ast,2,2,2,2,2), \quad \tau\sim(\ast,1^3,21,21,21,21)
\end{equation}

\noindent
with a placeholder \(\ast\) for the \textit{polarization} \((\varkappa_1,\tau_1)\).
The stabilization remains constant along the entire coclass tree,
whereas the polarization varies from branch to branch
with strictly increasing \(\tau_1\).
Therefore, it suffices to list the polarization.
The graph theoretic structure of the tree is determined by
the \textit{nuclear rank} \(\nu\) and the descendant numbers \(N_1/C_1\) with step size \(s=1\),
respectively \(N_2/C_2\) with step size \(s=2\), when \(\nu=2\).
The latter phenomenon, which is called \textit{bifurcation} to bigger coclass,
occurs two times in the irregular pre-period of length \(4\).
Starting with the first primitive period, also of length \(4\),
the value of the nuclear rank \(\nu=1\) remains settled,
and no further information is required in order to describe the entire infinite tree.
The \(p\)-multiplicator rank \(\mu\) coincides with the \textit{relation rank} \(d_2\)
in the Shafarevich theorem
\cite{Sh}.


\renewcommand{\arraystretch}{1.1}

\begin{table}[ht]
\caption{Structure and periodicity of the coclass tree \(\mathcal{T}^2\langle 3\,125,4\rangle\)}
\label{tbl:RootOrd3125Id4}
\begin{center}
\begin{tabular}{|rl|rl|r|cc|cc|lc|}
\hline
\multicolumn{2}{|c|}{id} & \multicolumn{2}{c|}{ord}        & cl     & \(\varkappa_1\) & \(\tau_1\) & \(\nu\) & \(\mu\) & \(N_1/C_1\) & \(N_2/C_2\) \\
\hline
     \(4\) &             &        \(3\,125\) &    \(=5^5\) &  \(3\) &           \(0\) &    \(1^3\) &   \(1\) &   \(3\) & \(7/7\)     &             \\
   \(564\) &  \(-\#1;1\) &       \(15\,625\) &    \(=5^6\) &  \(4\) &           \(0\) &    \(1^4\) &   \(2\) &   \(4\) & \(30/5\)    & \(75/75\)   \\
   \(888\) &  \(-\#1;4\) &       \(78\,125\) &    \(=5^7\) &  \(5\) &           \(0\) &   \(21^3\) &   \(1\) &   \(4\) & \(12/8\)    &             \\
   \(889\) &  \(-\#1;5\) &       \(78\,125\) &    \(=5^7\) &  \(5\) &           \(1\) &   \(21^3\) &   \(0\) &   \(3\) &             &             \\ 
   \(890\) &  \(-\#1;6\) &       \(78\,125\) &    \(=5^7\) &  \(5\) &           \(1\) &   \(21^3\) &   \(0\) &   \(3\) &             &             \\ 
           &  \(-\#1;6\) &      \(390\,625\) &    \(=5^8\) &  \(6\) &           \(0\) & \(2^21^2\) &   \(2\) &   \(5\) & \(38/7\)    & \(75/15\)   \\
\hline
           & \(-\#1;11\) &   \(1\,953\,125\) &    \(=5^9\) &  \(7\) &           \(0\) &   \(2^31\) &   \(1\) &   \(5\) & \(13/3\)    &             \\
           &  \(-\#1;6\) &   \(9\,765\,625\) & \(=5^{10}\) &  \(8\) &           \(0\) &    \(2^4\) &   \(1\) &   \(5\) & \(29/7\)    &             \\
           &  \(-\#1;1\) &  \(48\,828\,125\) & \(=5^{11}\) &  \(9\) &           \(0\) &   \(32^3\) &   \(1\) &   \(5\) & \(13/3\)    &             \\
           &  \(-\#1;6\) & \(244\,140\,625\) & \(=5^{12}\) & \(10\) &           \(0\) & \(3^22^2\) &   \(1\) &   \(5\) & \(37/7\)    &             \\

\hline
\end{tabular}
\end{center}
\end{table}


\renewcommand{\arraystretch}{1.1}

\begin{table}[ht]
\caption{Exemplary root path to \(S\simeq\mathrm{Gal}(\mathrm{F}_5^\infty(K)/K)\) with \(\mathrm{lo}=14\)}
\label{tbl:RootOrd15625Id564Lo14}
\begin{center}
\begin{tabular}{|rl|rl|r|r|cc|cc|lc|}
\hline
\multicolumn{2}{|c|}{id} & \multicolumn{2}{c|}{ord}           & cl    & dl    & \(\varkappa_1\) & \(\tau_1\) & \(\nu\) & \(\mu\) & \(N_1/C_1\) & \(N_2/C_2\) \\
\hline
   \(564\) &             &          \(15\,625\) &    \(=5^6\) & \(4\) & \(2\) &           \(0\) &    \(1^4\) &   \(2\) &   \(4\) & \(30/5\)    &  \(75/75\)  \\
           & \(-\#2;27\) &         \(390\,625\) &    \(=5^8\) & \(5\) & \(3\) &           \(1\) &    \(1^5\) &   \(1\) &   \(3\) &  \(3/3\)    &             \\
           &  \(-\#1;3\) &      \(1\,953\,125\) &    \(=5^9\) & \(6\) & \(3\) &           \(1\) &    \(1^5\) &   \(2\) &   \(4\) & \(10/5\)    &   \(5/5\)   \\
           &  \(-\#2;1\) &     \(48\,828\,125\) & \(=5^{11}\) & \(7\) & \(3\) &           \(1\) &    \(1^5\) &   \(1\) &   \(3\) & \(13/13\)   &             \\
           &  \(-\#1;6\) &    \(244\,140\,625\) & \(=5^{12}\) & \(8\) & \(3\) &           \(1\) &    \(1^5\) &   \(2\) &   \(4\) &  \(6/0\)    &   \(5/0\)   \\
           &  \(-\#2;1\) & \(6\,103\,515\,625\) & \(=5^{14}\) & \(9\) & \(3\) &           \(1\) &    \(1^5\) &   \(0\) &   \(2\) &             &             \\
\hline
\end{tabular}
\end{center}
\end{table}



In Table
\ref{tbl:RootOrd15625Id564Lo14},
resp.
\ref{tbl:RootOrd15625Id564Lo21},
we give an exemplary, resp. exceptional, \textit{root path}
to one of the non-metabelian Schur \(\sigma\)-groups \(S\simeq\mathrm{Gal}(\mathrm{F}_5^\infty(K)/K)\)
of order \(5^{14}\), derived length \(3\), resp. order \(5^{21}\), derived length \(4\), and corresponding invariants.
In Table
\ref{tbl:RootOrd15625Id564Lo21},
there occurs a \textit{trifurcation}
at the vertex with order \(5^{13}\), nuclear rank \(\nu=3\), and descendant numbers \(N_3/C_3\).


\renewcommand{\arraystretch}{1.1}

\begin{table}[ht]
\caption{Exceptional root path to \(S\simeq\mathrm{Gal}(\mathrm{F}_5^\infty(K)/K)\) with \(\mathrm{lo}=21\)}
\label{tbl:RootOrd15625Id564Lo21}
\begin{center}
\begin{tabular}{|rl|rl|r|r|cc|cc|lcc|}
\hline
\multicolumn{2}{|c|}{id} & \multicolumn{2}{c|}{ord}                  & cl     & dl    & \(\varkappa_1\) & \(\tau_1\) & \(\nu\) & \(\mu\) & \(N_1/C_1\) & \(N_2/C_2\) & \(N_3/C_3\) \\
\hline
   \(564\) &             &                 \(15\,625\) &    \(=5^6\) &  \(4\) & \(2\) &           \(0\) &    \(1^4\) &   \(2\) &   \(4\) & \(30/5\)    &  \(75/75\)  &             \\
           & \(-\#2;24\) &                \(390\,625\) &    \(=5^8\) &  \(5\) & \(3\) &           \(1\) &    \(1^5\) &   \(1\) &   \(3\) &  \(3/3\)    &             &             \\
           &  \(-\#1;3\) &             \(1\,953\,125\) &    \(=5^9\) &  \(6\) & \(3\) &           \(1\) &    \(1^5\) &   \(2\) &   \(4\) & \(30/5\)    &  \(25/25\)  &             \\
           &  \(-\#2;1\) &            \(48\,828\,125\) & \(=5^{11}\) &  \(7\) & \(3\) &           \(1\) &    \(1^5\) &   \(2\) &   \(4\) & \(78/78\)   & \(313/313\) &             \\
           &\(-\#2;306\) &        \(1\,220\,703\,125\) & \(=5^{13}\) &  \(8\) & \(3\) &           \(1\) &    \(1^5\) &   \(3\) &   \(5\) & \(39/5\)    &  \(95/50\)  &  \(25/25\)  \\
           &  \(-\#3;1\) &      \(152\,587\,890\,625\) & \(=5^{16}\) &  \(9\) & \(3\) &           \(1\) &    \(1^5\) &   \(1\) &   \(3\) & \(13/13\)   &             &             \\
           & \(-\#1;11\) &      \(762\,939\,453\,125\) & \(=5^{17}\) & \(10\) & \(3\) &           \(1\) &    \(1^5\) &   \(2\) &   \(4\) & \(10/0\)    &   \(5/5\)   &             \\
           &  \(-\#2;1\) &  \(19\,073\,486\,328\,125\) & \(=5^{19}\) & \(11\) & \(3\) &           \(1\) &    \(1^5\) &   \(1\) &   \(3\) &  \(3/3\)    &             &             \\
           &  \(-\#1;3\) &  \(95\,367\,431\,640\,625\) & \(=5^{20}\) & \(12\) & \(3\) &           \(1\) &    \(1^5\) &   \(1\) &   \(3\) &  \(5/0\)    &             &             \\
           &  \(-\#1;1\) & \(476\,837\,158\,203\,125\) & \(=5^{21}\) & \(13\) & \(4\) &           \(1\) &    \(1^5\) &   \(0\) &   \(2\) &             &             &             \\
\hline
\end{tabular}
\end{center}
\end{table}

\newpage


\begin{figure}[ht]
\caption{\(5\)-class field tower groups of order \(5^{14}\) on the tree \(\mathcal{T}\langle 3\,125,4\rangle\)}
\label{fig:SchurSigmaLogOrd14}


\setlength{\unitlength}{0.7cm}
\begin{picture}(18,26.5)(-6,-25.5)

{\tiny

\put(-5,0.5){\makebox(0,0)[cb]{Order}}
\put(-5,0){\line(0,-1){24}}
\multiput(-5.1,0)(0,-2){13}{\line(1,0){0.2}}
\put(-5.2,0){\makebox(0,0)[rc]{\(25\)}}
\put(-4.8,0){\makebox(0,0)[lc]{\(5^2\)}}
\put(-5.2,-2){\makebox(0,0)[rc]{\(125\)}}
\put(-4.8,-2){\makebox(0,0)[lc]{\(5^3\)}}
\put(-5.2,-4){\makebox(0,0)[rc]{\(625\)}}
\put(-4.8,-4){\makebox(0,0)[lc]{\(5^4\)}}
\put(-5.2,-6){\makebox(0,0)[rc]{\(3\,125\)}}
\put(-4.8,-6){\makebox(0,0)[lc]{\(5^5\)}}
\put(-5.2,-8){\makebox(0,0)[rc]{\(15\,625\)}}
\put(-4.8,-8){\makebox(0,0)[lc]{\(5^6\)}}
\put(-5.2,-10){\makebox(0,0)[rc]{\(78\,125\)}}
\put(-4.8,-10){\makebox(0,0)[lc]{\(5^7\)}}
\put(-5.2,-12){\makebox(0,0)[rc]{\(390\,625\)}}
\put(-4.8,-12){\makebox(0,0)[lc]{\(5^8\)}}
\put(-5.2,-14){\makebox(0,0)[rc]{\(1\,953\,125\)}}
\put(-4.8,-14){\makebox(0,0)[lc]{\(5^9\)}}
\put(-5.2,-16){\makebox(0,0)[rc]{\(9\,765\,625\)}}
\put(-4.8,-16){\makebox(0,0)[lc]{\(5^{10}\)}}
\put(-5.2,-18){\makebox(0,0)[rc]{\(48\,828\,125\)}}
\put(-4.8,-18){\makebox(0,0)[lc]{\(5^{11}\)}}
\put(-5.2,-20){\makebox(0,0)[rc]{\(244\,140\,625\)}}
\put(-4.8,-20){\makebox(0,0)[lc]{\(5^{12}\)}}
\put(-5.2,-22){\makebox(0,0)[rc]{\(1\,220\,703\,125\)}}
\put(-4.8,-22){\makebox(0,0)[lc]{\(5^{13}\)}}
\put(-5.2,-24){\makebox(0,0)[rc]{\(6\,103\,515\,625\)}}
\put(-4.8,-24){\makebox(0,0)[lc]{\(5^{14}\)}}
\put(-5,-24){\vector(0,-1){2}}


\put(6,0.5){\makebox(0,0)[lc]{Legend:}}
\put(5.8,-0.7){\framebox(0.4,0.4){}}
\put(6,-1){\circle*{0.2}}
\put(6,-1.5){\circle{0.2}}
\put(5.9,-2.1){\framebox(0.2,0.2){}}
\put(5.9,-2.6){\framebox(0.2,0.2){}}
\put(6,-2.5){\makebox(0,0)[cc]{\(\times\)}}
\put(5.9,-3.1){\framebox(0.2,0.2){\vrule height 3pt width 3pt}}
\put(5.95,-3.55){\framebox(0.1,0.1){}}

\put(6.5,-0.5){\makebox(0,0)[lc]{abelian}}
\put(6.5,-1){\makebox(0,0)[lc]{metabelian with coclass 1}}
\put(6.5,-1.5){\makebox(0,0)[lc]{metabelian with coclass 2}}
\put(6.5,-2){\makebox(0,0)[lc]{non-metabelian with TKT L}}
\put(6.5,-2.5){\makebox(0,0)[lc]{exceptional non-metabelian with TKT L}}
\put(6.5,-3){\makebox(0,0)[lc]{non-metabelian Schur \(\sigma\)-group with TKT L}}
\put(6.5,-3.5){\makebox(0,0)[lc]{non-metabelian with TKT M}}


\put(0,-16.5){\makebox(0,0)[cc]{Fork path symbol:}}
\put(0,-18){\makebox(0,0)[cc]{
\(\overbrace{\mathrm{L}\binom{1}{\rightarrow}}^{\text{Leaf}}\quad
\overbrace{\mathrm{M}}^{\text{Fork}}\quad
\overbrace{\left\lbrace\binom{2}{\leftarrow}\mathrm{L}\binom{1}{\leftarrow}\mathrm{L}\right\rbrace^{2}\ }^{\text{Trunk}}\quad
\overbrace{\binom{2}{\leftarrow}\mathrm{L}}^{\text{Leaf}}\)
}}
\put(0,-20){\makebox(0,0)[cc]{Transfer kernel types (TKT):}}
\put(0,-20.5){\makebox(0,0)[cc]{L: \(\varkappa=(122222)\), M: \(\varkappa_0=(022222)\)}}

\put(-1.2,-0.2){\framebox(0.4,0.4){}}
\put(-0.9,0.2){\makebox(0,0)[lb]{\(\langle 2\rangle\)}}
\put(-1,0){\line(0,-1){2}}
\put(-1,-2){\circle*{0.2}}
\put(-0.9,-1.8){\makebox(0,0)[lb]{\(\langle 3\rangle\)}}

\put(-1,-2){\vector(0,-1){2}}
\put(-1.2,-3.3){\makebox(0,0)[rc]{infinite}}
\put(-1.2,-3.8){\makebox(0,0)[rc]{mainline}}
\put(-1.2,-4.3){\makebox(0,0)[rc]{\(\mathcal{T}^1\langle 25,2\rangle\)}}

\put(-0.5,-2.8){\makebox(0,0)[lb]{metabelian bifurcation}}
\put(-1,-2){\line(1,-4){1}}

\put(0,-6){\line(0,-1){4}}
\multiput(0,-6)(0,-2){3}{\circle{0.2}}
\put(0.1,-5.8){\makebox(0,0)[lb]{\(\langle 4\rangle\)}}
\put(-0.1,-7.7){\makebox(0,0)[rb]{fork}}
\put(0.1,-7.8){\makebox(0,0)[lb]{\(\langle 564\rangle\)}}
\put(0.1,-9.8){\makebox(0,0)[lb]{\(\langle 888\rangle\)}}
\put(0,-8){\line(-1,-2){1}}
\put(0,-8){\line(-1,-1){2}}
\multiput(-1,-10)(-1,0){2}{\circle{0.2}}
\put(-1,-10.2){\makebox(0,0)[ct]{\(\langle 889\rangle\)}}
\put(-2,-10.2){\makebox(0,0)[ct]{\(\langle 890\rangle\)}}

\put(-1.5,-10){\oval(2.5,1.2)}
\put(-3.8,-8.5){\makebox(0,0)[lc]{\(S/S^{\prime\prime}\simeq\)}}
\put(-4.8,-9){\makebox(0,0)[lc]{\(M=\mathrm{Gal}(\mathrm{F}_5^2(K)/K)\)}}

\put(0,-10){\vector(0,-1){2}}
\put(-0.2,-11.3){\makebox(0,0)[rc]{infinite}}
\put(-0.2,-11.8){\makebox(0,0)[rc]{mainline}}
\put(-0.2,-12.3){\makebox(0,0)[rc]{\(\mathcal{T}^2\langle 3125,4\rangle\)}}

\put(2.5,-8.8){\makebox(0,0)[lb]{first non-metabelian bifurcation}}
\put(0,-8){\line(1,-4){1}}
\put(0,-8){\line(3,-5){2.2}}
\put(0,-8){\line(3,-1){12}}

\multiput(0.9,-12.1)(0.5,0){5}{\framebox(0.2,0.2){}}
\put(0.5,-11.8){\makebox(0,0)[cb]{\(2;\)}}
\put(1,-11.8){\makebox(0,0)[cb]{\(30\)}}
\put(1.5,-11.8){\makebox(0,0)[cb]{\(29\)}}
\put(2,-11.8){\makebox(0,0)[cb]{\(28\)}}
\put(2,-12){\makebox(0,0)[cc]{\(\times\)}}
\put(2.5,-11.8){\makebox(0,0)[cb]{\(27\)}}
\put(3,-11.8){\makebox(0,0)[cb]{\(26\)}}
\put(4.25,-12.8){\makebox(0,0)[cc]{\(\overbrace{\qquad\qquad\qquad 8\text{ times similarly } \qquad\qquad\qquad}\)}}
\multiput(5.4,-12.1)(0.5,0){5}{\framebox(0.2,0.2){}}
\put(5.5,-11.8){\makebox(0,0)[cb]{\(25\)}}
\put(6,-11.8){\makebox(0,0)[cb]{\(24\)}}
\put(6,-12){\makebox(0,0)[cc]{\(\times\)}}
\put(6.5,-11.8){\makebox(0,0)[cb]{\(23\)}}
\put(7,-11.8){\makebox(0,0)[cb]{\(22\)}}
\put(7.5,-11.8){\makebox(0,0)[cb]{\(21\)}}
\multiput(9.95,-12.05)(0.5,0){5}{\framebox(0.1,0.1){}}
\put(10,-11.8){\makebox(0,0)[cb]{\(20\)}}
\put(10.5,-11.8){\makebox(0,0)[cb]{\(19\)}}
\put(11,-11.8){\makebox(0,0)[cb]{\(18\)}}
\put(11.5,-11.8){\makebox(0,0)[cb]{\(17\)}}
\put(12,-11.8){\makebox(0,0)[cb]{\(16\)}}

\multiput(1,-12)(4.5,0){2}{\line(0,-1){2}}
\multiput(1.5,-12)(4.5,0){1}{\line(0,-1){2}}
\multiput(2.5,-12)(4.5,0){2}{\line(0,-1){2}}
\multiput(6.5,-12)(4.5,0){1}{\line(0,-1){2}}
\multiput(3,-12)(4.5,0){2}{\line(0,-1){2}}
\multiput(2.4,-14.1)(4.5,0){1}{\framebox(0.2,0.2){}}
\multiput(2.4,-13.8)(4.5,0){1}{\makebox(0,0)[rb]{\(1;3\)}}

\put(4.5,-14.8){\makebox(0,0)[lb]{second non-metabelian bifurcation}}
\put(2.5,-14){\line(1,-1){4}}
\put(2.5,-14){\line(2,-1){8}}
\multiput(6.4,-18.1)(1,0){5}{\framebox(0.2,0.2){}}
\put(6.5,-17.8){\makebox(0,0)[cb]{\(2;1\)}}
\put(7.5,-17.8){\makebox(0,0)[cb]{\(2;2\)}}
\put(8.5,-18.8){\makebox(0,0)[cc]{\(\overbrace{\qquad 5\text{ times similarly } \qquad}\)}}
\put(8.5,-17.8){\makebox(0,0)[cb]{\(2;3\)}}
\put(9.5,-17.8){\makebox(0,0)[cb]{\(2;4\)}}
\put(10.5,-17.8){\makebox(0,0)[cb]{\(2;5\)}}

\multiput(6.5,-18)(1,0){5}{\line(0,-1){2}}
\put(6.4,-20.1){\framebox(0.2,0.2){}}
\put(6.4,-19.8){\makebox(0,0)[rb]{\(1;7\)}}

\put(7.5,-20.8){\makebox(0,0)[lb]{third non-metabelian bifurcation}}
\put(6.5,-20){\line(1,-4){1}}
\put(6.5,-20){\line(5,-4){5}}

\multiput(7.4,-24.1)(1,0){5}{\framebox(0.2,0.2){\vrule height 3pt width 3pt}}
\put(7.5,-24.2){\makebox(0,0)[ct]{\(2;1\)}}
\put(8.5,-24.2){\makebox(0,0)[ct]{\(2;2\)}}
\put(9.5,-24.2){\makebox(0,0)[ct]{\(2;3\)}}
\put(10.5,-24.2){\makebox(0,0)[ct]{\(2;4\)}}
\put(11.5,-24.2){\makebox(0,0)[ct]{\(2;5\)}}
\put(9.5,-25){\makebox(0,0)[cc]{\(5\) Schur \(\sigma\)-groups with \(\mathrm{dl}=3\)}}

\put(9.5,-24.4){\oval(5.5,1.5)}
\put(6.8,-23.5){\makebox(0,0)[rc]{\(S=\mathrm{Gal}(\mathrm{F}_5^\infty(K)/K)\)}}

}

\end{picture}

\end{figure}

\newpage


\newpage


\section{Conclusion}
\label{s:Conclusion}

\noindent
For the prime numbers \(p\in\lbrace 5,7\rbrace\),
we have demonstrated impressively
that descendant trees \(\mathcal{T}\) of finite \(p\)-groups with
edges of fixed step size \(s=1\) and vertices of fixed coclass \(\mathrm{cc}=2\)
and \(p\)-power order
reveal some \textit{general structural principles} of so-called \textit{coclass trees}:

\begin{itemize}
\item
The branches \(\mathcal{B}(1),\ldots,\mathcal{B}(p-1)\)
form an irregular \textit{pre-period}
with \(\frac{p-1}{2}\) \textit{bifurcations} to the next coclass \(\mathrm{cc}=3\)
at the roots of the branches \(\mathcal{B}(2),\mathcal{B}(4),\ldots,\mathcal{B}(p-1)\) with even counter.
\item
With respect to their isomophism as subgraphs of the tree,
the branches show a \textit{periodicity}
\cite{As1,Ma2015a,Ma2015c,Ma2018b,Ma2018c,Ne1,Ne2}
with length \(\ell=p-1\), that is,

\begin{equation}
\label{eqn:Periodicity}
(\forall\,i\ge p) \quad \mathcal{B}(i+\ell)\simeq\mathcal{B}(i).
\end{equation}

\item
In view of the \textit{realization} of the vertices \(V\in\mathcal{T}\)
as second \(p\)-class groups \(V\simeq\mathrm{Gal}(\mathrm{F}_p^2(K)/K)\)
of imaginary quadratic number fields \(K=\mathbb{Q}(\sqrt{d})\) with \(d<0\),
the terminal leaves of the branches \(\mathcal{B}(2),\mathcal{B}(4),\ldots\)
with even counters are \textit{admissible}
\cite{Ma2012,Ma2010,Ma2014,Ma2013}.
\end{itemize}



\section{Acknowledgement}
\label{s:Acknowledgement}

\noindent
The author acknowledges that his research was supported by the
Austrian Science Fund (FWF): P 26008-N25 and J0497-PHY,
and by the Research Executive Agency of the European Union (EUREA).
Required computational tools are described in
\cite{AHL,BCP,BCFS,HEO,MAGMA,Nm,Ob}.





\begin{thebibliography}{XX}
%
\bibitem{Ag}
M. Arrigoni,
\textit{On Schur \(\sigma\)-groups},
Math. Nachr.
\textbf{192}
(1998),
71--89.
%
\bibitem{Ar1}
E. Artin,
\textit{Beweis des allgemeinen Reziprozit\"atsgesetzes},
Abh. Math. Sem. Univ. Hamburg
\textbf{5}
(1927),
353--363,
DOI 10.1007/BF02952531.
%
\bibitem{Ar2}
E. Artin,
\textit{Idealklassen in Oberk\"orpern und allgemeines Reziprozit\"atsgesetz},
Abh. Math. Sem. Univ. Hamburg
\textbf{7}
(1929),
46--51,
DOI 10.1007/BF02941159.
%
\bibitem{AHL}
J. A. Ascione, G. Havas and C. R. Leedham-Green,
\textit{A computer aided classification of certain groups of prime power order},
Bull. Austral. Math. Soc.
\textbf{17}
(1977),
257--274, Corrigendum 317--319, Microfiche Supplement p. 320,
DOI 10.1017/s0004972700010467.
%
\bibitem{As1}
J. A. Ascione,
\textit{On \(3\)-groups of second maximal class},
Ph.D. Thesis,
Australian National University,
Canberra,
1979.
\bibitem{As2}
J. A. Ascione,
\textit{On \(3\)-groups of second maximal class},
Bull. Austral. Math. Soc.
\textbf{21}
(1980),
473--474.
%
\bibitem{Bg}
G. Bagnera,
\textit{La composizione dei gruppi finiti il cui grado \`e la quinta potenza di un numero primo},
Ann. di Mat.
(Ser. 3)
\textbf{1}
(1898),
137--228,
DOI 10.1007/bf02419191.
%
\bibitem{BaBu}
L. Bartholdi and M. R. Bush,
\textit{Maximal unramified \(3\)-extensions of imaginary quadratic fields and \(\mathrm{SL}_2\mathbb{Z}_3\)},
J. Number Theory
\textbf{124}
(2007),
159--166.
%
\bibitem{BEO1}
H. U. Besche, B. Eick, and E. A. O'Brien,
\textit{A millennium project: constructing small groups},
Int. J. Algebra Comput.
\textbf{12}
(2002),
623-644.
%
\bibitem{BEO2}
H. U. Besche, B. Eick and E. A. O'Brien,
\textit{The SmallGroups Library --- a Library of Groups of Small Order},
2005,
an accepted and refereed GAP package, available also in MAGMA.
%
\bibitem{BCP}
W. Bosma, J. Cannon, and C. Playoust,
\textit{The Magma algebra system. I. The user language},
J. Symbolic Comput.
\textbf{24}
(1997),
235--265.
%
\bibitem{BCFS}
W. Bosma, J. J. Cannon, C. Fieker, and A. Steels (eds.),
\textit{Handbook of Magma functions}
(Edition 2.25,
Sydney,
2020).
%
\bibitem{BBH}
N. Boston, M. R. Bush and F. Hajir,
\textit{Heuristics for \(p\)-class towers of imaginary quadratic fields},
Math. Ann.
\textbf{368}
(2017),
no. 1,
633--669,
DOI 10.1007/s00208-016-1449-3.
%
\bibitem{Br}
J. R. Brink,
\textit{The class field tower for imaginary quadratic number fields of type \((3,3)\)},
Dissertation,
Ohio State University,
1984.
%
\bibitem{BrGo}
J. R. Brink and R. Gold, 
\textit{Class field towers of imaginary quadratic fields}, 
manuscripta math.
\textbf{57}
(1987),
425--450.
%
\bibitem{BuMa}
M. R. Bush and D. C. Mayer,
\textit{\(3\)-class field towers of exact length \(3\)},
J. Number Theory,
\textbf{147}
(2015),
766--777,
DOI 10.1016/j.jnt.2014.08.010.
%
\bibitem{GNO}
G. Gamble, W. Nickel, and E. A. O'Brien,
\textit{ANU p-Quotient --- p-Quotient and p-Group Generation Algorithms},
2006,
an accepted GAP package, available also in MAGMA.
%
\bibitem{GAP}
GAP Developer Group,
\textit{GAP -- Groups, Algorithms, and Programming --- a System for Computational Discrete Algebra},
Version 4.11.0,
Aachen, Braunschweig, Fort Collins, St. Andrews,
2020,\\
\texttt{(http://www.gap-system.org)}.
%
\bibitem{HeSm}
F.-P. Heider und B. Schmithals,
\textit{Zur Kapitulation der Idealklassen in unverzweigten primzyklischen Erweiterungen},
J. Reine Angew. Math.
\textbf{336}
(1982),
1--25.
%
\bibitem{HEO}
D. F. Holt, B. Eick, and E. A. O'Brien,
\textit{Handbook of computational group theory},
Discrete mathematics and its applications,
Chapman and Hall/CRC Press,
Boca Raton,
2005.
%
\bibitem{KoVe}
H. Koch und B. B. Venkov,
\textit{\"Uber den \(p\)-Klassenk\"orperturm eines imagin\"ar-quadra\-tischen Zahlk\"orpers},
Ast\'erisque
\textbf{24--25}
(1975),
57--67.
%
\bibitem{MAGMA}
MAGMA Developer Group,
\textit{MAGMA Computational Algebra System},
Version 2.25-4,
Sydney,
2020,\\
\texttt{(http://magma.maths.usyd.edu.au)}.
%
\bibitem{Ma1991}
D. C. Mayer,
\textit{Principalization in complex \(S_3\)-fields},
Congressus Numerantium
\textbf{80}
(1991),
73--87.
(Proceedings of the Twentieth Manitoba Conference on Numerical Mathematics and Computing,
Univ. of Manitoba, Winnipeg, Canada, 1990.)
%
\bibitem{Ma2012}
D. C. Mayer,
\textit{The second \(p\)-class group of a number field},
Int. J. Number Theory
\textbf{8}
(2012),
no. 2,
471--505,
DOI 10.1142/S179304211250025X.
%
\bibitem{Ma2010}
D. C. Mayer,
\textit{Transfers of metabelian \(p\)-groups},
Monatsh. Math.
\textbf{166}
(2012),
no. 3--4,
467--495,
DOI 10.1007/s00605-010-0277-x.
%
\bibitem{Ma2014}
D. C. Mayer,
\textit{Principalization algorithm via class group structure},
J. Th\'eor. Nombres Bordeaux
\textbf{26}
(2014),
no. 2,
415--464,
DOI 10.5802/jtnb.874.
%
\bibitem{Ma2013}
D. C. Mayer,
\textit{The distribution of second \(p\)-class groups on coclass graphs},
J. Th\'eor. Nombres Bordeaux
\textbf{25}
(2013),
no. 2,
401--456,
DOI 10.5802/jtnb.842.
(27i\`emes Journ\'ees Arithm\'etiques,
Faculty of Mathematics and Informatics,
University of Vilnius,
Lithuania,
July 2011.)
%
\bibitem{Ma2015a}
D. C. Mayer,
\textit{Periodic bifurcations in descendant trees of finite \(p\)-groups},
Adv. Pure Math.
\textbf{5}
(2015),
no. 4,
162--195,
DOI 10.4236/apm.2015.54020,
Special Issue on Group Theory,
March 2015.
%
\bibitem{Ma2015b}
D. C. Mayer,
\textit{Index-\(p\) abelianization data of \(p\)-class tower groups},
Adv. Pure Math.
\textbf{5}
(2015)
no. 5,
286--313,
DOI 10.4236/apm.2015.55029,
Special Issue on Number Theory and Cryptography,
April 2015.
%
\bibitem{Ma2015c}
D. C. Mayer,
\textit{Periodic sequences of \(p\)-class tower groups},
J. Appl. Math. Phys.
\textbf{3}
(2015),
no. 7,
746--756,
DOI 10.4236/jamp.2015.37090.
%
\bibitem{Ma2015d}
D. C. Mayer,
\textit{New number fields with known \(p\)-class tower},
Tatra Mt. Math. Pub.,
\textbf{64}
(2015),
21--57,
DOI 10.1515/tmmp-2015-0040,
Special Issue on Number Theory and Cryptology \lq 15.
%
\bibitem{Ma2016a}
D. C. Mayer,
\textit{Artin transfer patterns on descendant trees of finite \(p\)-groups},
Adv. Pure Math.
\textbf{6}
(2016),
no. 2,
66 -- 104,
DOI 10.4236/apm.2016.62008,
Special Issue on Group Theory Research,
January 2016.
%
\bibitem{Ma2016b}
D. C. Mayer,
\textit{\(p\)-Capitulation over number fields with \(p\)-class rank two},
J. Appl. Math. Phys.
\textbf{4}
(2016),
no. 7,
1280--1293,
DOI 10.4236/jamp.2016.47135.
%
\bibitem{Ma2016c}
D. C. Mayer,
\textit{\(p\)-Capitulation over number fields with \(p\)-class rank two},
\(2\)nd International Conference on Groups and Algebras (ICGA) \(2016\),
Suzhou, China,
invited lecture delivered on July 26, 2016, \\
\texttt{http://www.algebra.at/ICGA2016Suzhou.pdf}.
%
\bibitem{Ma2016d}
D. C. Mayer,
\textit{Recent progress in determining \(p\)-class field towers},
Gulf J. Math.
\textbf{4}
(2016),
no. 4,
74--102,
ISSN 2309-4966.
%
\bibitem{Ma2016e}
D. C. Mayer,
\textit{Recent progress in determining \(p\)-class field towers},
\(1\)st International Colloquium of Algebra, Number Theory, Cryptography and Information Security (ANCI) \(2016\),
Facult\'e Polydisciplinaire de Taza,
Universit\'e Sidi Mohamed Ben Abdellah,
F\`es, Morocco,
invited keynote delivered on November 12, 2016, \\
\texttt{http://www.algebra.at/ANCI2016DCM.pdf}.
%
\bibitem{Ma2017}
D. C. Mayer,
\textit{Criteria for three-stage towers of \(p\)-class fields}.
Adv. Pure Math.
\textbf{7}
(2017),
no. 2,
135--179,
DOI 10.4236/apm.2017.72008,
Special Issue on Number Theory,
February 2017.
%
\bibitem{Ma2018}
D. C. Mayer,
\textit{Deep transfers of \(p\)-class tower groups},
\(3\)rd International Conference on Groups and Algebras (ICGA) \(2018\),
Engineering Information Institute (EngII),
Sanya, Hainan, China, invited keynote delivered on January 14, 2018,
\texttt{http://www.algebra.at/DCM@ICGA2018Sanya.pdf}.
%
\bibitem{Ma2018a}
D. C. Mayer,
\textit{Deep transfers of \(p\)-class tower groups},
J. Appl. Math. Phys.
\textbf{6}
(2018),
no. 1,
DOI 10.4236/jamp.2018.61005.
%
\bibitem{Ma2018b}
D. C. Mayer,
\textit{Modeling rooted in-trees by finite \(p\)-groups},
Chapter 5, pp. 85--113,
in the Open Access Book \textit{Graph Theory --- Advanced Algorithms and Applications},
Ed. B. Sirmacek,
InTech d.o.o., Rijeka, January 2018,
DOI 10.5772/intechopen.68703.
%
\bibitem{Ma2018c}
D. C. Mayer,
\textit{Co-periodicity isomorphisms between forests of finite \(p\)-groups}.
Adv. Pure Math.
\textbf{8}
(2018),
no. 2,
77--140,
DOI 10.4236/apm.2018.81006,
Special Issue on Group Theory Studies,
January 2018.
%
\bibitem{Ma2018d}
D. C. Mayer,
\textit{Annihilator ideals of two-generated metabelian \(p\)-groups},
J. Algebra Appl.
\textbf{17}
(2018),
no. 4,
DOI 10.1142/S0219498818500767. 
%
\bibitem{Ma2020}
D. C. Mayer,
\textit{Pattern recognition via Artin transfers applied to class field towers},
\(3\)rd International Conference on Mathematics and its Applications (ICMA) \(2020\),
Facult\'e des Sciences d' Ain Chock Casablanca (FSAC), Universit\'e Hassan II,
Casablanca, Morocco, invited keynote delivered on February 28, 2020, \\
\texttt{http://www.algebra.at/DCM@ICMA2020Casablanca.pdf}.
%
\bibitem{Ne1}
B. Nebelung,
\textit{Klassifikation metabelscher \(3\)-Gruppen mit Faktorkommutatorgruppe vom Typ \((3,3)\)
und Anwendung auf das Kapitulationsproblem},
Inauguraldissertation, Band 1,
Universit\"at zu K\"oln,
1989.
%
\bibitem{Ne2}
B. Nebelung,
\textit{Anhang zu Klassifikation metabelscher \(3\)-Gruppen mit Faktorkommutatorgruppe vom Typ \((3,3)\)
und Anwendung auf das Kapitulationsproblem},
Inauguraldissertation, Band 2,
Universit\"at zu K\"oln,
1989.
%
\bibitem{Nm}
M. F. Newman,
\textit{Determination of groups of prime-power order},
pp. 73--84,
in: Group Theory, Canberra, 1975,
Lecture Notes in Math.,
vol. 573,
Springer,
Berlin,
1977.
%
\bibitem{Ob}
E. A. O'Brien, 
\textit{The \(p\)-group generation algorithm}, 
J. Symbolic Comput.
\textbf{9}
(1990),
677--698.
%
\bibitem{SoTa}
A. Scholz und O. Taussky,
\textit{Die Hauptideale der kubischen Klassenk\"orper imagin\"ar quadratischer Zahlk\"orper:
ihre rechnerische Bestimmung und ihr Einflu\ss\ auf den Klassenk\"orperturm},
J. Reine Angew. Math.
\textbf{171}
(1934),
19--41.
%
\bibitem{Sh}
I. R. Shafarevich,
\textit{Extensions with prescribed ramification points} (Russian),
Publ. Math., Inst. Hautes \'Etudes Sci. 
\textbf{18}
(1964),
71--95.
(English transl. by J. W. S. Cassels in
Amer. Math. Soc. Transl.,
II. Ser.,
\textbf{59}
(1966),
128--149.)

%
\end{thebibliography}
\end{document}